\documentclass[a4paper,11pt,reqno]{amsart}

\usepackage{dsfont,amssymb,color}

\usepackage{hyperref}

\numberwithin{equation}{section}

\newtheorem{theorem}{Theorem}[section]
\newtheorem{corollary}[theorem]{Corollary}
\newtheorem{lemma}[theorem]{Lemma}
\newtheorem{proposition}[theorem]{Proposition}
\theoremstyle{remark}
\newtheorem{remark}[theorem]{Remark}

\theoremstyle{definition}

\newcommand\bp{\begin{proof}}
\newcommand\ep{\end{proof}}

\newcommand\un{\mathds 1}

\newcommand{\C}{{\mathbb C}}
\newcommand{\R}{{\mathbb R}}
\newcommand{\Z}{{\mathbb Z}}

\newcommand{\F}{{\mathcal F}}

\newcommand{\B}{{\mathcal B}}

\newcommand\K{\mathcal K}
\newcommand\OO{\mathcal O}
\newcommand\RR{\mathcal R}

\newcommand\Sch{\mathcal S}

\newcommand\g{{\mathfrak g}}

\newcommand{\bj}{{\mathcal J}}

\newcommand{\vf}{\varphi}

\DeclareMathOperator\Dhat{{\hat\Delta}}
\DeclareMathOperator\Ad{\operatorname{Ad}}
\DeclareMathOperator{\arcsinh}{arcsinh}
\DeclareMathOperator{\Jac}{\mathrm{Jac}}

\begin{document}

\title{On deformations of C$^*$-algebras by actions of K\"{a}hlerian Lie groups}

\author[P. Bieliavsky]{Pierre Bieliavsky}

\email{Pierre.Bieliavsky@uclouvain.be}

\address{Institut de Recherche en Math\'ematique et Physique, Universit\'e Catholique de Louvain, Chemin du Cyclotron, 2, 1348 Louvain-la-Neuve, Belgium}

\author[V. Gayral]{Victor Gayral}

\email{victor.gayral@univ-reims.fr}

\address{Laboratoire de Math\'ematiques, Universit\'e de Reims Champagne-Ardenne,
Moulin de la Housse - BP 1039,
51687 Reims, France}

\author[S. Neshveyev]{Sergey Neshveyev}

\email{sergeyn@math.uio.no}

\address{Department of Mathematics, University of Oslo,
P.O. Box 1053 Blindern, NO-0316 Oslo, Norway}

\thanks{The research leading to these results has received funding from the European Research Council
under the European Union's Seventh Framework Programme (FP/2007-2013) / ERC Grant Agreement no. 307663%--NCGQG
}

\author[L. Tuset]{Lars Tuset}

\email{Lars.Tuset@hioa.no}

\address{Department of Computer Science, Oslo and Akershus University College of Applied Sciences,
P.O. Box 4 St. Olavs plass, NO-0130 Oslo, Norway}

\date{August 31, 2015; new version September 19, 2019\\
The previous version of the paper relied on unitarity of the dual cocycles we consider. As it turned out, these cocycles are only coisometric. Since the paper is already published, the present version corrects only Appendix~\ref{ap:A} and adds Appendix~\ref{ap:B} explaining why the main result remains true.}

\begin{abstract}
We show that two approaches to equivariant strict deformation quantization of C$^*$-algebras by actions of negatively curved K\"ahlerian Lie groups, one based on oscillatory integrals and the other on quantizations maps defined by dual $2$-cocycles, are equivalent.
\end{abstract}

\maketitle

\bigskip

\tableofcontents

\section*{Introduction}

The aim of this note is to establish equivalence of two approaches to equivariant strict deformation quantization of C$^*$-algebras equipped with actions of negatively curved K\"ahlerian Lie groups. The first approach is motivated by Rieffel's theory of deformation quantization for actions of $\R^d$~\cite{Ri1} and is based on the formalism of oscillatory integrals extended to these groups. This approach has recently been developed by the first two authors~\cite{BG} as a culmination of the program initiated in~\cite{Bi}. The second approach departs from the general theory of deformations of C$^*$-algebras by actions of locally compact quantum groups and dual measurable cocycles, developed by the third and fourth authors~\cite{NTdef2}. This theory, in turn, was motivated by Kasprzak's work~\cite{Kas} on deformation quantization for actions of abelian groups. It is known by now that for $\R^d$ the approaches of Rieffel and Kasprzak are equivalent~\cite{BNS,N}, but all available proofs rely crucially on commutativity of the group~$\R^d$. In particular, an important feature of deformations by actions of abelian groups is that the deformed algebras are equipped with actions of the same groups, while for non-abelian groups the symmetries of the deformed algebras should rather be quantum groups. This feature will be studied in detail in a subsequent
publication. Furthermore, non-unimodularity of the groups we consider pose an additional difficulty, in that the $*$-structures on dense subalgebras of the deformed C$^*$-algebras obtained by our two deformation procedures become incompatible. Our main result is that nevertheless the C$^*$-algebras are still canonically isomorphic. This gives, in our opinion, a sound justification of both deformation procedures. Our result also provides new tools for studying the deformed quantum groups, by combining the operator algebraic techniques suggested from the approach in~\cite{NTdef2} with the fine harmonic analysis of~\cite{BG} which, in particular, allows control at the smooth level too. This will be utilized in subsequent publications. We would also like to stress that most arguments are quite general, so there is reason to believe that once the results in~\cite{BG} are extended to a larger class of Lie groups, it should not take much effort to show compatibility with~\cite{NTdef2}.

\bigskip

\section{Preliminaries}
%\subsection{Generalities on K\"ahlerian Lie groups}
\label{GKLG}

From the seminal work \cite{PS} of Pyatetskii-Shapiro on bounded homogeneous
(not necessarily symmetric)
domains of $\C^n$ it is known that any K\"ahlerian Lie group with negative sectional curvature
(negatively curved, for short) can be written as an iterated semi-direct product
\begin{equation}
\label{KLG}
\big(\big(\dots\big(G_n\ltimes G_{n-1}\big)\ltimes\dots\big)\ltimes G_2\big)
\ltimes G_1
\end{equation}
of elementary blocks $G_j$ isomorphic to the Iwasawa factors $AN_j$ of the
simple Lie groups $SU(1,n_j)=KAN_j$.
Such blocks are called elementary K\"ahlerian Lie groups.
Hence, an elementary K\"ahlerian Lie group $G=AN$ is a solvable non-unimodular real Lie group of
dimension $2d+2$ with Lie algebra~$\g$ having a basis $H$, $\{X_j\}^{2d}_{j=1}$, $E$ satisfying the relations
$$
[H,E]=2E,\ \ [H,X_j]=X_j,\ \ [E,X_j]=0,\ \ [X_i,X_j]=(\delta_{i+d,j}-\delta_{i,j+d})E.
$$
The exponential map $\g\to G$ is a global diffeomorphism, and we will mainly be working in the global coordinate system given by the diffeomorphism
\begin{align}
\label{chart1}
\R\times\R^{2d}\times\R\ni(a,v,t)\mapsto \exp\big\{aH\big\}\exp\Big\{\sum^{2d}_{j=1}v_jX_j+tE\Big\}\in G.
\end{align}
The group law then takes the form
$$
(a,v,t)(a',v',t')=\big(a+a',e^{-a'}v+v',e^{-2a'}t+t'+\tfrac{1}{2}e^{-a'}\omega_0(v,v')\big),
$$
where $\omega_0(v,v')=\sum^d_{i=1}(v_iv'_{i+d}-v_{i+d}v'_i)$ is the standard symplectic form on $\R^{2d}$.
In this coordinate system the Lebesgue measure  on $\R^{2d+2}$ defines a (left) Haar measure $dg$ on $G$
with modular function
$$
\Delta_G(a,v,t)=e^{-(2d+2)a}.
$$
Our convention (opposite to the one in~\cite{BG}) for the modular function is such that the equality
$$
\int_Gf(gh)\,dg=\Delta_G(h)^{-1}\int_Gf(g)\,dg
$$
holds.

\smallskip

Let now $G$ be an arbitrary negatively curved K\"ahlerian Lie group with Pyatetskii-Shapiro
decomposition~\eqref{KLG}.
An important feature of Pyatetskii-Shapiro's theory is that the extension homomorphisms
at each step takes values in ${\rm Sp}(\R^{2d_j})$ if, as a manifold, $G_j=\R\times \R^{2d_j}\times \R$.
This implies in particular that under the global parametrization of $g\in G$ by $g=g_1\dots g_n$
with $g_i\in G_i$, the product of the Haar measures of the groups $G_i$ defines a Haar measure
on $G$.

Unless otherwise specified, the $L^p$-spaces
on $G$ will always be considered with respect to the Haar measure.
We denote by $\lambda$ and $\RR$ the left and right regular representations, and by
$\rho$ the unitarization of $\RR$:
\begin{equation}
\label{left-right-right}
(\lambda_gf)(g'):=f(g^{-1}g')\,,\qquad (\RR_gf)(g'):=f(g'g)
\,,\qquad \rho_g:=\Delta_G^{1/2}(g)\,\RR_g.
\end{equation}
By $\widetilde X$ and $\underline X$ we mean the left-invariant and right-invariant
vector fields on $G$ associated to the elements $X$ and $-X$ of $\g$, so
\begin{equation}
\label{invVF}
\widetilde X:=\frac d{dt}\Big|_{t=0}\,\RR_{e^{tX}}\,,\qquad  \underline X:=\frac d{dt}\Big|_{t=0}\,
\lambda_{e^{tX}}\;.
\end{equation}
We also extend this notation to the whole universal enveloping Lie algebra $\mathcal U(\g)$.

An important function space on $G$, denoted by~$\Sch(G)$, is the analogue of the Euclidean
Schwartz space, where the notion of regularity is associated to left-invariant vector fields, and
the decay is measured by the specific smooth function
$$
\mathfrak{d}_G:G\to\R^*_+\,,\quad g\mapsto
\sqrt{1\,+\,\|\Ad_{g}\|^2\,+\,\|\Ad_{g^{-1}}\|^2},
$$
where $\Ad$ denotes the adjoint action of $G$ on $\g$ and the norm is the operator norm on the finite
dimensional vector space $\g$ for any chosen Euclidean structure. We call $\mathfrak{d}_G$ the modular weight (not to be confused with the modular function
$\Delta_G$). By \cite[Lemma 2.4]{BG} we know that it is a sub-multiplicative
weight on $G$ (see also \cite[Definition 2.1]{BG}), which basically means that
it satisfies
$$
\Delta(\mathfrak{d}_G)\leq \mathfrak{d}_G\otimes\mathfrak{d}_G, \quad
|\widetilde{X}\mathfrak{d}_G|\leq C_{L,X}\, \mathfrak{d}_G,
\quad
|\underline{X}\mathfrak{d}_G|\leq C_{R,X} \,\mathfrak{d}_G,\qquad
\forall X\in \mathcal U(\g),
$$
for constants $C_{L,X},C_{R,X}$ depending only on $X\in\mathcal U(\g)$.
As shown in \cite[Lemma 3.27]{BG}, in the elementary case it is, up to scalar factors, bounded above and below by the function
$$
(a,v,t)\mapsto
\cosh a+\cosh2a+|v|(1+e^{2a}+\cosh a)+|t|(1+e^{2a}).
$$
The Schwartz space $\Sch(G)$
is defined as the Fr\'echet completion of $C^\infty_c(G)$ associated with the family of semi-norms
\begin{align}
\label{god-whaches-you}
f\mapsto\big\|\mathfrak{d}_G^n\,{\widetilde X}f\|_\infty
\end{align}
for all $n\in\mathbb N$ and $X\in\mathcal U(\g)$ (clearly, it suffices to consider only a basis in $\mathcal U(\g)$).

\begin{remark}
\label{sleep-on-the-sky}
Using among other things that $\mathfrak{d}_G^{-1}\in L^p(G)$ for $p>2d+1$, it is possible to show that one can use any other
$L^p$-norm ($1\le p\le\infty$) in the definition of the semi-norms \eqref{god-whaches-you} without modifying
the topology of $\Sch(G)$.
One  can also  replace the left-invariant vector
fields in \eqref{god-whaches-you} by their right-invariant counterparts \eqref{invVF}. This follows because the left-invariant vector
fields are linear combinations of right-invariant vector fields with coefficients given by smooth functions
which, together with their derivatives (in the sense of left- or right-invariant vector fields),
are bounded by a power of $\mathfrak d_G$, and vice versa.
\end{remark}

The Schwartz space $\Sch(G)$ is a nuclear  Fr\'echet algebra
stable under group inversion, and
the left and  right regular actions are strongly continuous. Obviously $C_c^\infty(G)\subset\Sch(G)\subset C_0(G)$
with continuous dense inclusions. When $G$ is elementary, the space $\Sch(G)$ is densely contained in the ordinary Schwartz space
$\Sch(\R^{2d+2})$ in the coordinates chart
\eqref{chart1}.

We will need the following result.
\begin{lemma} \label{lFvsS}
For any negatively curved K\"ahlerian Lie group $G$, the Schwartz space $\Sch(G)$ is a dense subspace of the Fourier algebra $A(G)$.
\end{lemma}
\bp
The left regular representation is strongly continuous on $\Sch(G)$. Since by
Remark \ref{sleep-on-the-sky} we may use right-invariant vector
fields instead of left-invariant ones, we see that $\Sch(G)$ is its own
subspace of smooth vectors for  $\lambda$. By the Dixmier-Malliavin Theorem
it follows that $\Sch(G)$ also coincides with its G{\aa}rding subspace, that is, we have
$\Sch(G)=\Sch(G)\ast \Sch(G)$ (finite sum of convolution products).
% is a finite sum of elements of the form  $\lambda(\eta)\vf=\eta\ast \vf$,
%where $\eta\in C^\infty_c(G)\subset \Sch(G)$ and $\vf\in\Sch(G)$.
This proves the lemma, since $A(G)=L^2(G)\ast L^2_\rho(G)$, where
$L^2_\rho(G)$ is the $L^2$-space on $G$ for the right Haar measure, and since $\Sch(G)$
is dense in both $L^2(G)$ and~$L^2_\rho(G)$.
%any element
%of $A(G)$ is of the form $f_1\ast \check{f_2}$, $f_1,  f_2\in L^2(G)$ and where $\check f(g):=f(g^{-1})$.
 \ep

Another important function space is the non-Abelian analogue of the Laurent Schwartz's space~$\B$:
\begin{equation}
\label{BCBG}
\B(G):=\big\{F\in C^\infty(G)\;:\; \big\|{\widetilde X}F\|_\infty<\infty,\ \forall X\in\mathcal U(\g)\big\}.
\end{equation}
When $G$ is elementary, we can equivalently define $\B(G)$ using the increasing sequence of norms
\begin{equation}
\label{eBCBG2}
\|F\|_{k}:=\max_{j+j_1+\dots j_{2d}+j'\leq k}\;\big\|\,{\widetilde H}^{j}\,
{\widetilde X_1}^{j_1}\dots {\widetilde X_{2d}}^{j_{2d}}\,{\widetilde E}^{j'}F\|_\infty.
\end{equation}

It is shown in \cite[Lemma 2.8]{BG} that $\B(G)$ is Fr\'echet. In fact,
it coincides with the space
of smooth vectors for the right regular action $\RR$ within $C_{ru}(G)$, the $C^*$-algebra of right-uniformly
continuous and bounded functions on $G$. (Our convention for the right uniform structure on a
group is the one that yields strong continuity for the right regular action.) However, as opposed
to the Schwartz space~$\Sch(G)$, one cannot use right-invariant vector fields
to topologize  $\B(G)$ and it is  not stable under the group inversion.

\smallskip

Let us finally say a few words about elementary K\"ahlerian Lie groups $G$. They are endowed with extra geometrical structures not shared by non-elementary ones.
Namely, they are also left $G$-equivariant symplectic symmetric spaces. By this
we mean that each $g\in G$ has a smooth involution $s_g:G\to G$
(the symmetry at $g$), having $g$ as a unique isolated fixed point, such that
$$
s_g\circ s_{g'} \circ s_g=s_{s_g(g')},
$$
together with a symplectic $2$-form $\omega$ on $G$ that is invariant
$s_g^\star\,\omega=\omega$
under the symmetries, and such that $\lambda$ acts by symplectomorphisms on $(G,\omega)$ in a covariant fashion
$$
\lambda_g\circ s_{g'}=s_{g^{-1}g'}\circ\lambda_g
$$
with respect to the symmetries. In the coordinates \eqref{chart1} the symmetries are given by
\begin{align*}
 s_{(a,v,t)}(a',v',t')&:=
\big(2a-a',2v\cosh(a-a')-v',2t\cosh(2a-2a')-t'+\omega_0(v,v')\sinh(a-a')\big),
\end{align*}
while the invariant symplectic form is given by $\omega:=2da\wedge dt\,+\,\omega_0$.
As a symplectic symmetric space, $G$ has a unique midpoint map, that is,
a smooth map ${\rm mid}:G\times G \to G$
such that $s_{{\rm mid}(g,g')}(g)=g'$ for all $g,g'\in G$.
Moreover, the medial triangle  map
\begin{align}
\label{mie-de-pain}
\Phi_G: G^3\to G^3\,,\quad (g_1,g_2,g_3)\mapsto \big({\rm mid}(g_1,g_2),{\rm mid}(g_2,g_3),
{\rm mid}(g_3,g_1)\big),
\end{align}
is a global diffeomorphism invariant under the diagonal left action of $G$.

We also mention the decomposition  $G=Q\ltimes P$, which reflects the existence
of a global (real) polarization on the symplectic manifold $(G,\omega)$, where
\begin{align}
\label{QP}
Q=\exp\Big\{\R H+\sum_{j=1}^{d} \R X_j\Big\}\quad\mbox{and}\quad
P=\exp\Big\{\sum_{j=d+1}^{2d} \R X_j+\R E\Big\}.
\end{align}
The group  $Q$ is non-unimodular and solvable, while $P$ is Abelian.
%It will also be useful to consider the second global coordinate system associated to this
%decomposition:
%\begin{align}
%\label{chart2}
%\R\times\R^d\times\R^d\times\R\ni(a,x,y,t)
%\mapsto \exp\Big\{aH+ \sum^{d}_{j=1}x_jX_j\Big\}\exp\Big\{\sum^{d}_{j=1}y_jX_{j+d}+tE\Big\}\in G.
%\end{align}

\bigskip

\section{Deformations of function algebras}

In this section we fix an elementary K\"ahlerian Lie group $G$. We aim to compare the two deformations of function algebras of $G$ studied in~\cite{BG} and~\cite{NTdef2}.

\subsection{Deformations of \texorpdfstring{$C_0(G)$}{C0(G)}}

For a fixed parameter $\theta\in\R^*$ consider the two-point kernel on $G$ defined by
\begin{align}
\label{KTheta}
K_\theta(g_1,g_2)=\frac{4}{(\pi\theta)^{2d+2}}\,A(g_1,g_2)\exp\Big\{{\tfrac{2i}{\theta}S(g_1,g_2)}\Big\},
\end{align}
where, with $\Phi_G$ the medial triangle map given in \eqref{mie-de-pain}, we define
\begin{align*}
S(g_1,g_2):={\rm Area}\left(\Phi_G^{-1}(e,g_1,g_2)\right)\,,\quad
A(g_1,g_2):=\mbox{\rm Jac}_{\Phi_G^{-1}}^{1/2}(e,g_1,g_2).
\end{align*}
Here ${\rm Area}(g_1,g_2,g_3)$ is the symplectic area of any surface in $G$ admitting
an oriented geodesic triangle $T(g_1,g_2,g_3)$ as boundary. (This is unambiguously defined
since as a manifold $G$ has trivial de Rham cohomology in degree two.)
In the coordinates \eqref{chart1}, with $g_j=(a_j,v_j,t_j)$, we have
$$
A(g_1,g_2)=\big(\cosh(a_1)\cosh(a_2)\cosh(a_1-a_2)\big)^d\big(\cosh(2a_1)\cosh(2a_2)
\cosh(2a_1-2a_2)\big)^{1/2},
$$
$$
S(g_1,g_2)=\sinh(2a_1)t_2-\sinh(2a_2)t_1+\cosh(a_1)\cosh(a_2)\omega_0(v_1,v_2).
$$
It is sometimes useful to consider the associated three point kernel
$$
K_\theta^3(g_1,g_2,g_3):=K_\theta(g_1^{-1}g_2,g_1^{-1}g_3),
$$
which of course is invariant under the diagonal left action of $G$. But since the functions
$$
{\rm Area}\left(\Phi_G^{-1}(g_1,g_2,g_3)\right)\quad \mbox{and}\quad
\mbox{\rm Jac}_{\Phi_G^{-1}}(g_1,g_2,g_3)
$$
are also invariant under the diagonal left action of $G$ as well as under cyclic permutations,
we see that $K^3_G$ is also invariant under cyclic permutations.
At the level of the two point kernel this implies
\begin{align}
\label{weird-inv}
K_\theta(g^{-1},g^{-1}h)=K_\theta(h,g).
\end{align}
In passing we record another important symmetry property
\begin{align}
\label{not-weird-inv}
\overline{K_{\theta}(g,h)}=K_{-\theta}(g,h)=K_\theta(h,g).
\end{align}

By \cite[Proposition 4.10]{BG} the formula below endows $\Sch(G)$ with a new involutive
and associative Fr\'echet
algebra structure (the involution is still complex conjugation
and the topology is unaltered):
\begin{align}
\label{SP}
f_1\star_{\theta} f_2=\int_{G\times G}K_\theta(g_1,g_2)\,\RR_{g_1}(f_1)\,\RR_{g_2}(f_2)\,dg_1\,dg_2.
\end{align}
Property \eqref{not-weird-inv} entails
\begin{equation} \label{esym}
f_1\star_{-\theta} f_2=f_2\star_{\theta} f_1.
\end{equation}

Up to a nontrivial unitary transformation of $L^2(G)$ that commutes with complex conjugation,
this deformed product is, in the chart \eqref{chart1}, the usual Moyal product  on $\R^{2d+2}$.
More precisely, by \cite[Theorem 6.43 and Lemma 7.10]{BG} there is a $G$-equivariant quantization map ${\rm Op}_{G,\theta}$,
denoted by $\Omega_{\theta,{\bf m}_0}$ in~\cite{BG}, which defines a unitary operator from $L^2(G)$ to the Hilbert algebra of Hilbert-Schmidt operators on $L^2(Q)$.
Here $Q$ is the subgroup of $G$ entering the decomposition $G=Q\ltimes P$ given in~\eqref{QP}, and the space $L^2(Q)$ is equipped with an irreducible unitary representation $U_\theta$ of $G$, see \cite[Section~7.2]{BG}.
% of $G$ defined by
%$$
%U_\theta(a, v, t)\psi(a_0, n_0)
%= \exp\left\{\frac{i}{\theta}\Big(e^{2(a-a_0)}t + \omega_0\big(\frac{1}{2}e^{a-a_0}n - n_0, e^{a-a_0}m\big)\Big)\right\}\psi(a_0 - a, n_0 - e^{a-a_0}n),
%$$
%where we write the coordinates of $Q$ as $(a_0,n_0)\in\R\times\R^d$.
Then
$$f_1\star_{\theta} f_2={\rm Op}_{G,\theta}^{-1}\big({\rm Op}_{G,\theta}(f_1){\rm Op}_{G,\theta}(f_2)\big),
$$
and the required unitary transformation of $L^2(G)$
is given by ${\rm Op}_{W,\theta}^{-1}\circ{\rm Op}_G$, where ${\rm Op}_{W,\theta}$ denotes the Weyl quantization map. This transformation is $T_{\theta,0}^{-1}$ in the notation of \cite{BG}, and its explicit form is given in~\cite[Equation (62)]{BG}.

It follows that the deformed
product \eqref{SP} extends to the space $L^2(G)$, which then becomes a Hilbert algebra
isomorphic to the algebra of Hilbert-Schmidt operators on the separable Hilbert space $L^2(Q)$. In particular, we have a representation~$\pi_\theta$ of $(\Sch(G),\star_\theta)$ on $L^2(G)$ given by
$$
\pi_\theta(f_1)f_2=f_1\star_\theta f_2\ \ \text{for}\ \ f_1,f_2\in\Sch(G).
$$
The operators $\pi_\theta(f)$ are bounded, with $\|\pi_\theta(f)\|\leq \|f\|_2$, and satisfy
$\pi_\theta(f)^*=\pi_\theta(\bar f)$. Of course, this also implies that
 the $C^*$-algebra generated by $\pi_\theta(\Sch(G))$ is isomorphic to
the algebra of compact operators on $L^2(Q)$. This $C^*$-algebra is a deformation of $C_0(G)$, which we
coin $C_0(G)_\theta$.  Be aware that  $L^2(G)\subset C_0(G)_\theta$
but $C_0(G)\not\subset C_0(G)_\theta$ (or more precisely, $\pi_\theta$ extends to $L^2(G)$ but not to~$C_0(G)$).
This definition of $C_0(G)_\theta$ is slightly different from the one in
\cite[Proposition 8.26]{BG}, but it is
equivalent to it, in that we use the representation $\pi_\theta$ on $L^2(G)$ instead of the quasi-equivalent irreducible
representation on $L^2(Q)$ employed in \cite{BG}.

\medskip

Starting from the product $\star_\theta$ there is another natural construction of a C$^*$-algebra
deforming~$C_0(G)$, see~\cite{NTdef2}.
Let $W^*(G)$ be the von Neumann algebra generated by the image of the left regular representation
$\lambda$ on $L^2(G)$. As usual the Fourier algebra $A(G)$
 is  identified  with the pre\-dual~$W^*(G)_*$ of~$W^*(G)$ using the pairing
$(f,\lambda_g)=f(g)$.
%The important observation is that the space $A(G)$ is closed under the (suitably
%extended) product~$\star_\theta$.

Recall~\cite[Section 5.1]{NTdef2} that the kernel $K_\theta$
defines a dual unitary $2$-cocycle $\Omega_\theta$ on $G$, initially  defined as the quadratic form
on $\Sch(G\times G)$ given by
$$
\big(\Omega_\theta\xi,\zeta\big):=\int_{G\times G}\overline{K_\theta(g_1,g_2)}\,
\big((\lambda_{g_1^{-1}}\otimes\lambda_{g_2^{-1}})\xi,\zeta\big)
\,dg_1\,dg_2.
$$
(Our scalar products are linear in the first variable.) As it was not proven in~\cite{NTdef2} that this form indeed defines a unitary operator on $L^2(G\times G)$, for the reader's convenience we supply a possible argument in Appendix~\ref{ap:A}. Thus, $\Omega_\theta$ is a unitary element in $W^*(G)\bar\otimes W^*(G)$ satisfying the
cocycle identity
$$
(\Omega_\theta\otimes1)(\Dhat\otimes\iota)(\Omega_\theta)=(1\otimes\Omega_\theta)(\iota\otimes\Dhat)
(\Omega_\theta),
$$
where $\Dhat\colon W^*(G)\to W^*(G)\bar\otimes W^*(G)$ is the comultiplication defined by
$\Dhat(\lambda_g)=\lambda_g\otimes\lambda_g$.

Since $\Sch(G)=\Sch(G)*\Sch(G)\subset A(G)$ by Lemma~\ref{lFvsS}, we have
$$
f_1\star_\theta f_2=(f_1\otimes f_2)(\Dhat(\cdot)\Omega_\theta^*)\ \ \text{for}\ \ f_1,f_2\in \Sch(G)\subset A(G).
$$
This identity can be used to extend $\star_\theta$ to the whole space $A(G)$, but we are not going to do this and will always work
with the dense subspace $\Sch(G)$ of $A(G)$.

Consider now the multiplicative unitary $\hat W\in W^*(G)\bar\otimes L^\infty(G)$ of the dual quantum group $\hat G$, so
\begin{equation}
\label{MU}
(\hat W\xi)(g,h)=\xi(hg,h)=(\lambda_h^{-1}\xi(\cdot,h))(g)\ \ \text{for}\ \ \xi\in L^2(G\times G).
\end{equation}
According to \cite[Sections~2.1 \&~3.1]{NTdef2} we can define a representation $\pi_{\Omega_\theta}$ of $(\Sch(G),\star_\theta)$ on $L^2(G)$ by
$$
\pi_{\Omega_\theta}(f)=(f\otimes\iota)(\hat W\Omega_\theta^*).
$$
The norm closure of $\pi_{\Omega_\theta}(\Sch(G))$ becomes a C$^*$-algebra, which we denote by $C^*_r(\hat G;\Omega_\theta)$.

\smallskip

In order to compare the algebras $C_0(G)_\theta$ and $C^*_r(\hat G;\Omega_\theta)$, consider the respective modular conjugations $J$ and~$\hat J$ of the group~$G$ and the dual quantum group $\hat G$, so
$$
(J\xi)(g)=\overline{\xi(g)},\ \ (\hat J\xi)(g)=\Delta_G^{-1/2}(g)\,\overline{\xi(g^{-1})}.
$$
Then, as already observed in \cite[Sections~4.1 \&~5.2]{NTdef2}, it follows from our definitions that
\begin{equation}\label{ekeyid0}
\pi_{\Omega_\theta}(f_1)\check f_2=(f_1\star_\theta f_2)\check{}\ \ \text{for}\ \ f_1,f_2\in \Sch(G),
\end{equation}
where $\check f(g)=f(g^{-1})$. Consider the involutive unitary given by the product
of the two modular conjugations
$$
\bj:=J\hat J=\hat JJ.
$$
Then \eqref{ekeyid0} implies that
\begin{equation} \label{ekeyid}
\pi_{\Omega_\theta}(f)=\bj\Delta_G^{-1/2}\pi_\theta(f)\Delta_G^{1/2}\bj\ \ \text{for}\ \ f\in \Sch(G).
\end{equation}
Here we view the modular function $\Delta_G$ as the (unbounded) operator of multiplication by $\Delta_G$
on~$L^2(G)$.

We are going to show that  $\Delta_G$  coincides on $\Sch(G)$ with the adjoint action of a
$\star_\theta$-multiplier of~$\Sch(G)$, for which
 we need  to introduce the following pseudo-differential operator on $G$:
$$
T_\theta:=\Big(
\big(1-\pi^2\theta^2\partial_t^2\big)^{1/2}+i\pi\theta\partial_t\Big)^{d+1}.
$$
We observe that $T_\theta$ commutes with left translations (as $\partial_t$ coincides, in the chart \eqref{chart1}, with the left-invariant vector field $\widetilde E$ associated to the element $E\in\g$), that it preserves the space $\Sch(G)$, and that $T_\theta^{-1}=T_{-\theta}$.
\begin{lemma}
\label{l:key}
Let $\alpha\in\C$.
The maps $f\mapsto  L_{\star_\theta}(\Delta_G^\alpha)f:=\Delta_G^\alpha \star_\theta f$ and $f\mapsto
R_{\star_\theta}(\Delta_G^\alpha)f:= f\star_\theta \Delta_G^\alpha  $ define invertible operators on $\Sch(G)$ which factorize as
$$
L_{\star_\theta}(\Delta_G^\alpha)=\Delta_G^\alpha\circ T_\theta^\alpha= T_\theta^\alpha\circ \Delta_G^\alpha\quad\mbox{and}\quad
R_{\star_\theta}(\Delta_G^\alpha)=\Delta_G^\alpha\circ T_{-\theta}^\alpha= T_{-\theta}^\alpha\circ \Delta_G^\alpha.
 $$
\end{lemma}

Here the expressions $\Delta_G^\alpha \star_\theta f$ and $f\star_\theta \Delta_G^\alpha  $ are defined by interpreting \eqref{SP} as an oscillatory integral, as explained in \cite[Chapter~4]{BG}. This requires  $\Delta^\alpha_G$ to be a tempered weight, which indeed
follows from the discussion in \cite{BG} just before Definition 2.6 and from Lemma 2.21 there. Furthermore, the operators $L_{\star_\theta}(\Delta_G^\alpha)$ and $R_{\star_\theta}(\Delta_G^\alpha)$ are continuous on $\Sch(G)$ by \cite[Proposition~4.10]{BG}.

\bp[Proof of Lemma~\ref{l:key}]
We only need to prove the decomposition
$L_{\star_\theta}(\Delta_G^\alpha)=\Delta_G^\alpha\circ T_\theta^\alpha$. Indeed, since
$\Delta_G$ only depends on the variable $a$, while $T_\theta$
is a continuous function of $i\partial_t$, the maps $\Delta_G^\alpha$ and  $T_\theta^\alpha$
commute. Hence $L_{\star_\theta}(\Delta_G^\alpha)= T_\theta^\alpha\circ \Delta_G^\alpha$. Next, the decomposition $L_{\star_\theta}(\Delta_G^\alpha)=\Delta_G^\alpha\circ T_\theta^\alpha$ also yields invertibility of $L_{\star_\theta}(\Delta_G^\alpha)$, since
$T_\theta^{-1}=T_{-\theta}$. Last, the relations for $R_{\star_\theta}(\Delta_G^\alpha)$
also follow, since $L_{\star_{-\theta}}(\Delta_G^\alpha)=R_{\star_\theta}(\Delta_G^\alpha)$, which, in turn, is
a consequence of $K_{-\theta}(g_1,g_2)=K_\theta(g_2,g_1)$.

To prove the factorization $L_{\star_\theta}(\Delta_G^\alpha)=\Delta_G^\alpha\circ T_\theta^\alpha$,
note first that the left-invariance of the deformed product $\star_\theta$ implies that the operator
$ \Delta_G^{-\alpha}\circ L_{\star_\theta}( \Delta_G^\alpha )$
commutes with left translations, whence it is of the form~$\RR(S)$ for a
distribution $S\in C^\infty_c(G)'\simeq C^\infty_c(\R^{2d+2})'$. To determine explicitly this
distribution, we proceed with  formal
computations which, however, can easily be made rigorous. With $g=(a,v,t)$, we have
\begin{align*}
\Delta_G^{-\alpha}(g)\big(\Delta_G^\alpha\star_\theta f\big)(g)&= \Delta_G^{-\alpha}(g)
\int K_\theta(g_1,g_2)\,\Delta_G^\alpha(gg_1)\,f(gg_2)\,dg_1\,dg_2\\&
=\int K_\theta(g_1,g_2)\,\Delta_G^\alpha(g_1)\,f_2(gg_2)\,dg_1\,dg_2\\
&=\frac{4}{(\pi\theta)^{2d+2}}\int A(a_1,a_2) e^{\tfrac{2i}\theta(\sinh(2a_1)t_2-\sinh(2a_2)t_1+\cosh(a_1)\cosh(a_2)\omega_0(v_1,v_2))}e^{\alpha(2d+2)a_1}
\\&\times
 f\big(a+a_2,e^{-a_2}v+v_2,e^{-2a_2}t+t_2+\tfrac{1}{2}e^{-a_2}\omega_0(v,v_2)\big)\,
da_1\,dv_1\,dt_1\,da_2\,dv_2\,dt_2\\
&=\frac{4}{\pi^2\theta^2}\int \frac{\cosh(a_1-a_2)^d}{\big(\cosh(a_1)\cosh(a_2)\big)^d}
\big(\cosh(2a_1)\cosh(2a_2)
\cosh(2a_1-2a_2)\big)^{1/2}\\
&\times e^{\alpha(2d+2)a_1}\,e^{\tfrac{2i}\theta(\sinh(2a_1)t_2-\sinh(2a_2)t_1)}
 f\big(a+a_2,e^{-a_2}v,e^{-2a_2}t+t_2\big)\,
da_1\,dt_1\,da_2\,dt_2\\
&=\frac{2}{\pi\theta}\int \cosh(2a_1)\,e^{\alpha(2d+2)a_1}\,e^{\tfrac{2i}\theta\sinh(2a_1)t_2}
 f\big(a,v,t+t_2\big)\,
da_1\,dt_2\\&
=\int \big(\pi\theta a_1+(1+(\pi\theta a_1)^2)^{1/2}\big)^{{\alpha(d+1)}
}\,e^{2i\pi a_1(t_2-t)}
 f(a,v,t_2)\,
da_1\,dt_2,
\end{align*}
which concludes the proof.
 \end{proof}

\begin{remark}
From the above Lemma it easily follows that $\Delta_G^\alpha\star_\theta\Delta_G^\beta=\Delta_G^{\alpha
+\beta}$ for all $\alpha,\beta\in\C$.
%In particular, we have for $t\in\R$:
%$$
% L_{\star_\theta}(\Delta_G^{it})=L_{\star_\theta}(\Delta_G)^{it}\quad\mbox{and}\quad
% R_{\star_\theta}(\Delta_G^{it})=R_{\star_\theta}(\Delta_G)^{it}.
%$$
\end{remark}

\begin{proposition} \label{pkey}
We have $C^*_r(\hat G;\Omega_\theta)=\bj C_0(G)_\theta\bj$.
\end{proposition}

\bp
From Lemma \ref{l:key} we deduce the following equalities for operators on $\Sch(G)$:
$$
\Delta_G^\alpha=T_\theta^{-\alpha}\circ L_{\star_\theta}(\Delta_G^\alpha)= L_{\star_\theta}(\Delta_G^\alpha)
\circ T_\theta^{-\alpha}\quad\mbox{and}\quad
\Delta_G^\alpha=T_{-\theta}^{-\alpha}\circ R_{\star_\theta}(\Delta_G^\alpha)=
 R_{\star_\theta}(\Delta_G^\alpha)
\circ T_{-\theta}^{-\alpha}.
$$
Since $ T_{-\theta}= T_{\theta}^{-1}$, we  then get
\begin{align*}
\Delta_G^{2\alpha}=L_{\star_\theta}(\Delta_G^\alpha)\circ R_{\star_\theta}(\Delta_G^\alpha)
=R_{\star_\theta}(\Delta_G^\alpha)\circ L_{\star_\theta}(\Delta_G^\alpha).
\end{align*}
With $\Delta_G^\alpha$ viewed as a
densely defined operator on $L^2(G)$ preserving its domain $\Sch(G)$,
the  relation above immediately implies
\begin{equation}\label{singing-on-the-train}
  \Delta_G^{-1/2}\pi_\theta(f)\Delta_G^{1/2}
=\pi_\theta(\Delta_G^{-1/4}\star_\theta f\star_\theta\Delta_G^{1/4}),
\end{equation}
as operators on $\Sch(G)$.
From this it follows that the map sending the operator $\pi_\theta(f)$ to the closure of the operator
$\Delta_G^{-1/2}\pi_\theta(f)\Delta_G^{1/2}$ defines an automorphism of the Fr\'echet space $\Sch(G)$,
 identified with $\pi_\theta(\Sch(G))$. The result is then an immediate consequence of
identity \eqref{ekeyid}.
\ep

Define an action $\beta$ of~$G$ on
$C^*_r(\hat G;\Omega_\theta)$ by $\beta_g=\Ad\rho_g$. Recall, see~\cite[Section~2.4]{NTdef2}, that the cocycle~$\Omega_\theta$ is called regular if $C^*_r(\hat G;\Omega_\theta)\rtimes_\beta G$ is isomorphic to the
algebra of compact operators on some Hilbert space. The condition of regularity plays an important role in the theory developed in~\cite{NTdef2}. As follows from the recent work of Baaj and Crespo~\cite{BC}, for general locally compact quantum groups this condition is equivalent to regularity of $G$, so in our case it is satisfied for any dual cocycle. This result is proved using the theory of quantum groupoids. For the cocycle $\Omega_\theta$, here is a more direct proof.

\begin{corollary}\label{cregular}
The dual cocycle $\Omega_\theta$ is regular.
\end{corollary}

\bp The isomorphism $C^*_r(\hat G;\Omega_\theta)\cong C_0(G)_\theta$, $x\mapsto \bj x\bj$,
intertwines the action $\beta$ with the action $\Ad\lambda_g$ on $C_0(G)_\theta$.
Specializing~\cite[Corollary~8.49]{BG} to  $A=\C$, we know that the crossed product
$C_0(G)_\theta\rtimes_{\Ad\lambda} G$ is Morita equivalent
to $\C$, hence it is isomorphic to the algebra of compact operators on a Hilbert space.
\ep

Recall also that the cocycle $\Omega_\theta$ is called continuous
if $\Omega_\theta\in M(C_r^*(G)\otimes C_r^*(G))$.

\begin{proposition}\label{pcontinuous}
The dual cocycle $\Omega_\theta$ is continuous.
\end{proposition}
\bp
We claim that both $\Omega_\theta$ and $\Omega_\theta^*$ preserve the space
$\Sch(G\times G)$. This is true by a minor extension of \cite[Lemma 2.49]{BG}, where
instead of the map $\mathcal R\otimes \mathcal R$ from \cite[Lemma 2.42]{BG}, one
considers the maps
\begin{align*}
&C^\infty(G\times G)\to
C^\infty\big(G\times G,C^\infty(G\times G)\big)\;,\\
&f \mapsto
\Big[(x,y)\in G\times G\mapsto (\lambda_x\otimes \lambda_y) (f)
:=\big[(g,h)\in G\times G \mapsto f(x^{-1}g,y^{-1}h)\big]\Big],\\
&f \mapsto
\Big[(x,y)\in G\times G\mapsto (\lambda_{x^{-1}}\otimes \lambda_{y^{-1}}) (f)
:=\big[(g,h)\in G\times G \mapsto f(xg,yh)\big]\Big].
\end{align*}

The proof then follows from a minor modification of \cite[Proposition 4.5]{NTdef2} using
the fact that $\Sch(G)$ is dense in $L^2(G)$.
\ep

\subsection{Oscillatory integrals and quantization maps}\label{s:qm}

Both deformation procedures work for a larger class of functions than $\Sch(G)$. Let us start by explaining the approach in \cite{BG}.

The product $\star_\theta$
extends to the space $\B(G)$ defined by \eqref{BCBG}, by replacing the ordinary integrals
appearing in~\eqref{SP} with oscillatory ones, see Chapters 2-4 in \cite{BG}.
Moreover, then $(\Sch(G),\star_\theta)$ is an ideal of $(\B(G),\star_\theta)$ and the representation $\pi_\theta$ extends (necessarily uniquely) to $(\B(G),\star_\theta)$. Namely,
$$
\pi_\theta(f)\xi=f\star_\theta\xi\ \ \text{for}\ \ f\in\B(G),\ \xi\in\Sch(G).
$$
By \cite[Theorems 8.20
\& 8.33]{BG} we have
\begin{align}
\label{to-get-sunlight}
\|\pi_\theta(f)\|\leq C\|f\|_K\ \ \text{for}\ \ f\in\B(G),
\end{align}
where
$C>0$ and $\|\cdot\|_K$ is one of the semi-norms~\eqref{eBCBG2} of $\B(G)$ with $K\in\mathbb N$ depending only on~${\rm dim}(G)$.

The representation $\pi_\theta$ can actually be described without using oscillatory integrals. In order to see this, we need
an important property of the product $\star_\theta$ called strong traciality, which means that under the
integral, deformed and pointwise products coincide. As the simple proof of this was omitted in~\cite{BG}, we include it here.
\begin{lemma}
\label{trace-of-foxes}
For $f_1,f_2\in\Sch(G)$, we have
$$
\int_G (f_1\star_\theta f_2)(g)\,dg=\int_G f_1(g)\, f_2(g)\,dg.
$$
\end{lemma}
\bp For any $f_j\in\Sch(G)$, $j=1,2,3$, we have
$$
(f_1\star_\theta f_2,f_3)=(f_2,\bar f_1\star_{\theta} f_3).
$$
By \cite[Proposition 5.19]{BG} (or rather its proof), any bounded approximate unit for the
commutative algebra $\Sch(G)$ is also a bounded approximate unit for the
non-commutative algebra $(\Sch(G),\star_\theta)$. Hence, letting $f_3$ run through
such an approximate unit, establishes the required equality.
\ep

We now give a description of $\pi_\theta$ using ordinary integrals.

\begin{lemma}\label{l:os-free}
For any $f\in\B(G)$ and $\xi,\zeta\in \Sch(G)$, we have
\begin{equation*} \label{eq:rep-pi}
(\pi_\theta(f)\xi,\zeta)=\int_G
f(g)\,(\xi\star_\theta\bar\zeta)(g)\,dg.
\end{equation*}
\end{lemma}

\bp
For $f\in\Sch(G)$ the result follows from Lemma~\ref{trace-of-foxes}. To extend it to all $f\in\B(G)$, we choose a uniformly bounded sequence $\{f_n\}_n$ in $\Sch(G)$ such that $f_n\to f$ in $\B^{\mathfrak d_G}(G)$, that is,
$$
\max_{j+j_1+\dots j_{2d}+j'\leq k}\;\big\|\mathfrak{d}_G^{-1}\,{\widetilde H}^{j}\,
{\widetilde X_1}^{j_1}\dots {\widetilde X_{2d}}^{j_{2d}}\,{\widetilde E}^{j'}\,(f-f_n)\|_\infty\to 0 \ \ \text{for all}\
k\in\mathbb N,
$$
which is possible by \cite[Lemma~2.8 (viii)]{BG}. Then $f_n\star_\theta \xi\to f\star_\theta \xi$ in $\Sch(G)$ by \cite[Theorem~4.9]{BG}, and the lemma follows by the dominated convergence theorem.
\ep

Let us now turn to the approach in \cite{NTdef2}.
Let $\hat W_{\Omega_\theta}$ be the multiplicative unitary of the locally compact quantum group $\hat G_{\Omega_\theta}$ defined as the von Neumann algebra $L^\infty(\hat G_{\Omega_\theta})=W^*(G)$ with the coproduct $\Dhat_{\Omega_\theta}=\Omega_\theta\Dhat(\cdot)\Omega_\theta^*$ and
invariant weights as defined by De Commer \cite{DC}. Using this unitary we can define `quantization
maps'
$$
T_\nu\colon L^\infty(G)\to \B(L^2(G)),\ \
\ \ f\mapsto(\iota\otimes\nu)\big(\hat W_{\Omega_\theta}\,\Omega_\theta\,(f\otimes1)\,
\Omega_\theta^*\,\hat W_{\Omega_\theta}^*\big),
$$
for $\nu\in \K(L^2(G))^*=\B(L^2(G))_*$. Here we identify a function $f\in L^\infty(G)$ with the operator of
multiplication by $f$ on $L^2(G)$. From now on we write $\K$ for the algebra of compact operators
$\K(L^2(G))$. It is shown in \cite[Lemma 3.2]{NTdef2} that
$$
C^*_r(\hat G;\Omega_\theta)=[T_\nu(f): f\in C_0(G),\ \nu\in \K^*],
$$
where the brackets $[\ ]$ denote closed linear span.

In order to exhibit the quantization maps more explicitly, recall that by~\cite[Proposition 5.4]{DC}
the multiplicative unitary~$\hat W_{\Omega_\theta}$ is given, with $\hat W$ as in \eqref{MU}, by
\begin{equation*}\label{eDC1}
\hat W_{\Omega_\theta}=(\tilde J\otimes \hat J)\Omega_\theta\hat W^*(J\otimes\hat J)\Omega_\theta^*.
\end{equation*}
The involution $\tilde J$ here (denoted by $J_N$ in \cite{DC}) is defined as follows. Consider the von Neumann algebra
$W^*(\hat G;\Omega_\theta)$ generated by~$C^*_r(\hat G;\Omega_\theta)$. The action $\beta$ extends
 to this von Neumann algebra by the same formula as before, so $\beta_g=\Ad\rho_g$. This action is
 integrable and ergodic,
hence it defines a n.s.f.~weight $\tilde\varphi$ on $W^*(\hat G;\Omega_\theta)$~by
$$
\tilde\varphi(x)1=\int_G\beta_g(x)\,dg\ \ \text{for}\ \ x\in W^*(\hat G;\Omega_\theta)_+.
$$
The space $L^2(G)$ can be identified with the space of the GNS-representation defined by $\tilde\varphi$. Namely, letting as usual ${\mathfrak N}_{\tilde\varphi}=\{x\mid\tilde\varphi(x^*x)<\infty\}$,
we have an $L^2$-norm isometric map $\tilde\Lambda\colon{\mathfrak N}_{\tilde\varphi}\to L^2(G)$ uniquely determined by
$$
\tilde\Lambda((f\otimes\iota)(\hat W\Omega_\theta^*))=(f\otimes\iota)(\hat W)
$$
for $f\in A(G)$ such that the right hand side is in $L^2(G)$. In other words, $\tilde\Lambda$ is uniquely defined by
$$
\tilde\Lambda(\pi_{\Omega_\theta}(f))=\check f\ \ \text{for}\ \ f\in \Sch(G).
$$
Then $\tilde J$ is the corresponding modular conjugation. Denote the
associated modular operator by $\tilde\Delta$.

\begin{proposition}\label{pmodular}
For the modular conjugation we have $\tilde J=J$, while the modular operator $\tilde\Delta$ is the closure of the operator
$$
f\mapsto(\Delta_G^{-1}\star_\theta \check f\star_\theta \Delta_G)\check{},\ \ f\in\Sch(G).
$$
In particular, the modular group of $\tilde\varphi$ is given by $\sigma^{\tilde\varphi}_t=
\Ad \Delta_G^{2it}$.
\end{proposition}

\bp It is more convenient to work with the von Neumann algebra $L^\infty(G)_\theta$ generated by
$C_0(G)_\theta$ and equipped with the action $\Ad\lambda_g$.
Since $C_0(G)_\theta$ is isomorphic to $\mathcal K(L^2(Q))$, $L^\infty(G)_\theta$
is isomorphic to $\B(L^2(Q))$.
Denote by $\tilde\varphi_\theta$ the corresponding weight on $L^\infty(G)_\theta$, so
$$
\tilde\varphi_\theta(x)1=\int_G(\Ad\lambda_g)(x)\,dg\ \ \text{for}\ \ x\in L^\infty(G)_{\theta,+}.
$$
Then $\tilde\varphi_\theta\circ(\Ad\lambda_g)=\Delta_G(g)^{-1}\tilde\varphi_\theta$.

On the other hand,
let $\psi_\theta$ be the operator trace on $\B(L^2(Q))$  transported to $L^\infty(G)_\theta$.
We want to  express $\tilde\varphi_\theta$ in terms of $\psi_\theta$.
For this, we first use the fact that ${\rm Op}_{G,\theta}$ (see the paragraph following
equation \eqref{esym}) is a unitary operator
from $L^2(G)$ to the Hilbert space of Hilbert-Schmidt operators on $L^2(Q)$, so for any $f\in\Sch(G)$ we have
$$
\psi_\theta(\pi_\theta(\bar f\star_\theta f))={\rm Tr}\big(\big|{\rm Op}_{G,\theta}(f)\big|^2\big)=\int_G|f(g)|^2\,dg.
$$
Denote by $\Delta_\theta$ the image of $\Delta_G$ under $\pi_\theta$. More precisely, we define
$\Delta_\theta$ as the generator of the one-parameter unitary group
$\{L_{\star_\theta}(\Delta^{it})\}_{t\in\R}=
\{L_{\star_\theta}
(\Delta)^{it}\}_{t\in\R}$. Therefore $\Delta_\theta$ is a positive unbounded operator affiliated with $L^\infty(G)_\theta$. Note that it easily follows from Lemma~\ref{l:key} that for any $f\in \Sch(G)$ the map $\alpha\mapsto\Delta_G^\alpha\star_\theta f\in L^2(G)$ is analytic. Hence $\Sch(G)$ is a core for $\Delta_\theta$. Consider now the weight
$$
\tilde\psi_\theta=\psi_\theta(\Delta_\theta^{-1/2}\cdot\Delta_\theta^{-1/2}).
$$
Since the product $\star_\theta$ is invariant under left translations, we have $(\Ad\lambda_g)(\Delta_\theta)=\Delta_G(g)^{-1}\Delta_\theta$. It follows that $\tilde\psi_\theta\circ (\Ad\lambda_g)=\Delta_G(g)^{-1}\tilde\psi_\theta$. This already implies that $\tilde\varphi_\theta= c\tilde\psi_\theta$ for some $c>0$. Indeed, as both $\tilde\varphi_\theta$ and $\tilde\psi_\theta$ are scaled the same way under the action $\Ad\lambda$, Connes' Radon-Nikodym cocycle $[D\tilde\varphi_\theta : D\tilde\psi_\theta]_t$ is $\Ad\lambda$-invariant. Since the action $\Ad\lambda$ on $L^\infty(G)_\theta$ is ergodic, the cocycle must be scalar-valued, and this implies the claim. We will see soon that $c=1$.

We can now identify $L^2(G)$ with the space of the GNS-representation defined by $\tilde\varphi_\theta$ using the isometric map $\tilde\Lambda_\theta\colon{\mathfrak N}_{\tilde\varphi_\theta}\to L^2(G)$ uniquely determined by
$$
\tilde\Lambda_\theta(\pi_\theta(f))=c^{1/2}f\star_\theta\Delta_G^{-1/2}\ \ \text{for}\ \ f\in\Sch(G).
$$
The corresponding modular operator $\tilde\Delta_\theta$ satisfies $\tilde\Delta_\theta^{it}\tilde\Lambda_\theta(x)=\tilde\Lambda_\theta(\Delta_\theta^{-it}x\Delta_\theta^{it})$, hence it is the closure of the operator
$$
f\mapsto \Delta_G^{-1}\star_\theta f\star_\theta\Delta_G, \ \ f\in\Sch(G).
$$
Since the modular conjugation satisfies $\tilde J_\theta\tilde\Delta^{1/2}_\theta\tilde\Lambda_\theta(x)=\tilde\Lambda_\theta(x^*)$ for $x\in{\mathfrak N}_{\tilde\varphi_\theta}\cap {\mathfrak N}_{\tilde\varphi_\theta}^*$, we also have $\tilde J_\theta=J$.

\smallskip

Let us return to $W^*(\hat G;\Omega_\theta)=\bj L^\infty(G)_\theta\bj$. Since by identities \eqref{ekeyid} and \eqref{singing-on-the-train} we have
\begin{equation}\label{e:keyid2}
\pi_{\Omega_\theta}(f)=\bj\pi_\theta(\Delta_G^{-1/4}\star_\theta f\star_\theta\Delta_G^{1/4})\bj\ \ \text{for}\ \ f\in \Sch(G),
\end{equation}
we can identify $L^2(G)$ with the space of the GNS-representation defined by $\tilde\varphi$ using the map
$$
\pi_{\Omega_\theta}(f)\mapsto \bj\tilde\Lambda_\theta(\pi_\theta(\Delta_G^{-1/4}\star_\theta f\star_\theta\Delta_G^{1/4}))=c^{1/2}(\Delta_G^{1/2}(\Delta_G^{-1/4}\star_\theta f\star_\theta\Delta_G^{-1/4}))\check{}
=c^{1/2}\check f.
$$
Comparing this with $\tilde\Lambda$ we conclude that $c=1$. It follows then that
$$
\tilde\Delta=\bj\tilde\Delta_\theta\bj,\ \ \tilde J=\bj\tilde J_\theta
\bj.
$$
Therefore $\tilde J=J$ and $\tilde\Delta$ is the closure of the operator
$$
\check f\mapsto (\Delta_G^{-1}\star_\theta f\star_\theta\Delta_G)\check{}, \ \ f\in\Sch(G).
$$
Finally, the statement about the modular group follows from
$$
\Delta_\theta^{-it}\pi_\theta(f)\Delta_\theta^{it}=\pi_\theta(\Delta_G^{-it}\star_\theta f\star_\theta\Delta_G^{it})=\Delta_G^{-2it}\pi_\theta(f)\Delta_G^{2it},
$$
as $\bj\Delta_G\bj=\Delta_G^{-1}$.
\ep

As a corollary we get
\begin{equation}\label{emult}
\hat W_{\Omega_\theta}\Omega_\theta=(J\otimes \hat J)\Omega_\theta\hat W^*(J\otimes\hat J).
\end{equation}

For $f\in L^1(G)$, we let as usual $\lambda(f)=\int_G f(g)\,\lambda(g)\,dg$.
\begin{lemma}
For $f_1,f_2\in\Sch(G)$, consider the function
$$
\eta_{f_1,f_2}=\Delta_G^{3/4}\star_\theta\bar f_2\star_\theta f_1\star_\theta\Delta_G^{1/4}\in\Sch(G).
$$
Then for the linear functional $\omega_{f_1,f_2}=(\cdot\,f_1,f_2)\in \K^*$ we have
$$
(\iota\otimes\omega_{f_1,f_2})(\hat W_{\Omega_\theta}\Omega_\theta)=\lambda(\check\eta_{f_1,f_2}).
$$
\end{lemma}

\bp Take $\xi\in \Sch(G)$. Using~\eqref{emult} we compute (the integrals below are
absolutely convergent):
\begin{align*}
(\hat W_{\Omega_\theta}\Omega_\theta(\xi\otimes f_1))(x,y)&=\Delta_G(y)^{-1/2}\overline{(\Omega_\theta\hat W^*(J\otimes\hat J)(\xi\otimes f_1))(x,y^{-1})}\\
&=\Delta_G(y)^{-1/2}\int K_\theta(g,h)\,\overline{(\hat W^*(J\otimes\hat J)(\xi\otimes f_1))(gx,hy^{-1})}\,dg\,dh\\
&=\Delta_G(y)^{-1/2}\int K_\theta(g,h)\,\overline{(J\xi\otimes \hat Jf_1)(yh^{-1}gx,hy^{-1})}\,dg\,dh\\
&=\int K_\theta(g,h)\,\Delta_G(h)^{-1/2}\,\xi(yh^{-1}gx)\,f_1(yh^{-1})\,dg\,dh.
\end{align*}
It follows that

\smallskip
$
\displaystyle\big((\iota\otimes\omega_{f_1,f_2})(\hat W_{\Omega_\theta}\Omega_\theta)\xi\big)(x)
$
\begin{align*}
&=\int K_\theta(g,h)\,\Delta_G(h)^{-1/2}\,(\lambda_{g^{-1}hy^{-1}}\xi)(x)f_1(yh^{-1})\,\overline{f_2(y)}\,dg\,dh\,dy\\
&=\int K_\theta(g,h)\,\Delta_G(h)^{-1/2}\,\Delta_G(y)^{-1}\,(\lambda_{g^{-1}hy}\xi)(x)\,f_1(y^{-1}h^{-1})\,\overline{f_2(y^{-1})}\,dg\,dh\,dy\\
&=\int K_\theta(g,h)\,\Delta_G(h)^{-1/2}\,\Delta_G(h^{-1}gy)^{-1}\,(\lambda_y\xi)(x)\,f_1(y^{-1}g^{-1})\,\overline{f_2(y^{-1}g^{-1}h)}\,dg\,dh\,dy\\
&=\int \check\eta(y)\,(\lambda_y\xi)(x)\,dy,
\end{align*}
where
\begin{align*}
\eta(y)&=\int K_\theta(g,h)\,(\Delta_G^{1/2}f_1)(yg^{-1})\,
(\Delta_G^{1/2}\bar f_2)(yg^{-1}h)\,dg\,dh\\
&=\int K_\theta(g^{-1},h)\,\Delta_G(y)\,(\Delta_G^{-1/2}f_1)(yg)\,
(\Delta_G^{1/2}\bar f_2)(ygh)\,dg\,dh\\
&=\int K_\theta(g^{-1},g^{-1}h)\,\Delta_G(y)\,(\Delta_G^{-1/2}f_1)(yg)\,
(\Delta_G^{1/2}\bar f_2)(yh)\,dg\,dh.
\end{align*}
Using \eqref{weird-inv} we get $\eta=\Delta_G((\Delta_G^{1/2}\bar f_2)\star_\theta(\Delta_G^{-1/2}f_1))$,
and from \eqref{singing-on-the-train} we deduce
$$
\eta
=\Delta_G^{3/4}\star_\theta\bar f_2\star_\theta f_1\star_\theta\Delta_G^{1/4}=\eta_{f_1,f_2},
$$
proving the lemma.
\ep

Similarly to $\lambda(f)$, for $f\in\Sch(G)$ we define an operator $\RR(f)=\int_G f(g)\,\RR_g\,dg$ acting on functions on $G$.

\begin{lemma}\label{l:I-have-no-money}
For any $f,f_1,f_2\in \Sch(G)$ we have
$$
T_{\omega_{f_1,f_2}}(\check f)=\pi_{\Omega_\theta}(\RR(\check\eta_{f_1,f_2})f).
$$
\end{lemma}

\bp By identity (3.1) in \cite{NTdef2} we have
\begin{equation*}
(\hat W_{\Omega_\theta}\Omega_\theta)_{23}\hat W_{12}(\hat W_{\Omega_\theta}\Omega_\theta)^*_{23}=(\hat W\Omega_\theta^*)_{12}(\hat W_{\Omega_\theta}\Omega_\theta)_{13}.
\end{equation*}
Applying $\iota\otimes\iota\otimes\omega_{f_1,f_2}$ to this, by the previous lemma we get
$$
(\iota\otimes T_{\omega_{f_1,f_2}})(\hat W)=\hat W\Omega_\theta^*(\lambda(\check\eta_{f_1,f_2})\otimes1).
$$
Applying now $f\otimes\iota$, we get the required identity, as $(f\otimes\iota)(\hat W)=\check f$ and the equality
$$
f(\cdot\,\lambda(\eta))=\RR(\eta)f
$$
holds in $A(G)$ for any $\eta\in L^1(G)$.
\ep

The maps $T_\nu\colon L^\infty(G)\to \B(L^2(G))$ are obviously ultraweakly continuous. On the other hand, for the representation $\pi_\theta\colon\B(G)\to  \B(L^2(G))$ we have the following result.

\begin{lemma}
\label{something-in-my-eye}
For any $\eta\in \Sch(G)$, the operator $\RR(\eta)$ maps $L^\infty(G)$ into $\B(G)$ and the map $L^\infty(G)\to \B(L^2(G))$, $f\mapsto \pi_\theta(\RR(\eta)f)$, is ultraweakly continuous.
\end{lemma}
\bp
Take $f\in L^\infty(G)$. Let $X\in\g$ and let $\widetilde X$ and $\underline X$ be the associated left-
and right-invariant vector fields defined in \eqref{invVF}. Then we find that
\begin{align*}
\widetilde X\big(\RR(\eta)f\big)(g)&=\frac d{dt}\Big|_{t=0}\big(\RR(\eta)f\big)(ge^{tX})
=\frac d{dt}\Big|_{t=0}\int_G \eta(g')\,f(ge^{tX}g')\,dg'\\
&=\frac d{dt}\Big|_{t=0}\int_G \eta(e^{-tX}g')\,f(gg')\,dg'=\int_G \big(\underline X\eta\big)(g')\,f(gg')\,dg'
=(\RR(\underline X\eta)f)(g),
\end{align*}
where we used the dominated convergence to exchange $t$-derivatives and integrals.
By induction we get a similar relation for every $X\in\mathcal U(\g)$ and thus finally arrive at the estimates
\begin{equation}\label{eq:I-have-no-names}
\|\widetilde X\big(R(\eta)f\big)\|_\infty\leq \|\underline X\eta\|_1\|f\|_\infty\,,
\quad \forall X\in\mathcal U(\g).
\end{equation}
Therefore $\RR(\eta)f\in\B(G)$.

Set $S_\eta(f):=\pi_\theta(R(\eta)f)$, so $S_\eta$ is a map $L^\infty(G)\to\B(L^2(G))$.
By inequalities \eqref{eq:I-have-no-names} and \eqref{to-get-sunlight} this map is bounded. Therefore in order to show that it is ultraweakly continuous it suffices to check that for any $\xi,\zeta\in\Sch(G)$ the linear functional $f\mapsto (S_\eta(f)\xi,\zeta)$ on $L^\infty(G)$ is ultraweakly continuous. By Lemma~\ref{l:os-free} we have
$$
(S_\eta(f)\xi,\zeta)=\int (\RR(\eta)f)(g)\,(\xi\star_\theta\bar\zeta)(g)\,dg
=\int f(gh)(\xi\star_\theta\bar\zeta)(g)\,\eta(h)\,dg\,dh.
$$
Since $(\xi\star_\theta\bar\zeta)\otimes\eta
\in\Sch(G)\otimes_{\mathrm{alg}}\Sch(G)\subset L^1(G\times G)$ and $\Delta\colon L^\infty(G)\to L^\infty(G\times G)$,
$\Delta(f)(g,h)=f(gh)$,
is ultraweakly continuous, we see that the linear functional $f\mapsto (S_\eta(f)\xi,\zeta)$ on $L^\infty(G)$ is indeed ultraweakly continuous.
\ep

For $\eta\in \Sch(G)$ put
$$
\eta^\theta=(\Delta_G^{-1/4}\star_\theta\eta\star_\theta\Delta_G^{1/4})\check{}.
$$
In particular, we have
$$
\eta^\theta_{f_1,f_2}=(\Delta_G^{3/4}\star_\theta\bar f_2\star_\theta f_1\star_\theta\Delta_G^{1/4})^\theta
=(\Delta_G^{1/2}\star_\theta\bar f_2\star_\theta f_1\star_\theta\Delta_G^{1/2})\check{}
=\Delta_G^{-1}(\bar f_2\star_\theta f_1)\check{}.
$$
We are now ready to describe the quantization maps in terms of $\pi_\theta$, which is the main result of this section.

\begin{proposition}\label{pquantmaps}
For any $f\in L^\infty(G)$ and $f_1,f_2\in \Sch(G)$ we have
$$
T_{\omega_{f_1,f_2}}(\check f)=\bj\pi_\theta\big(\RR(\eta^\theta_{f_1,f_2})f\big)\bj.
$$
\end{proposition}

\bp
Since both sides of the identity in the formulation are ultraweakly continuous in $f$, it suffices to check it for
$f\in \Sch(G)$. By Lemma~\ref{l:I-have-no-money} it is then enough to show that
$$
\pi_{\Omega_\theta}(\RR(\check\eta_{f_1,f_2})f)=\bj\pi_\theta(\RR(\eta^\theta_{f_1,f_2})f)\bj\ \ \text{for
all}\ \
f,f_1,f_2\in \Sch(G).
$$
Let us show that, a bit more generally,
$$
\pi_{\Omega_\theta}(\RR(\check\eta)f)=\bj\pi_\theta(\RR(\eta^\theta)f)\bj\ \ \text{for all}\ \ f,\eta\in\Sch(G).
$$
By identity \eqref{e:keyid2} the left-hand side equals
$$
\bj\pi_\theta(\Delta_G^{-1/4}\star_\theta (\RR(\check\eta)f)\star_\theta\Delta_G^{1/4})\bj.
$$
Therefore it remains to check
$$
\Delta_G^{-1/4}\star_\theta (\RR(\check\eta)f)\star_\theta\Delta_G^{1/4}=\RR(\eta^\theta)f.
$$

We have $\RR(\check\eta)f=\lambda(f)\eta$. Since $\star_\theta$ is invariant under left translations, we also have
$$
\Delta_G^{-1/4}\star_\theta (\lambda_g\eta)\star_\theta\Delta_G^{1/4}=\lambda_g(\Delta_G^{-1/4}\star_\theta\eta\star_\theta\Delta_G^{1/4}).
$$
Integrating with respect to the finite measure $f(g)dg$, the right-hand side becomes
$$
\lambda(f)(\Delta_G^{-1/4}\star_\theta\eta\star_\theta\Delta_G^{1/4})=\lambda(f){\eta^\theta}\,\check{}=
\RR(\eta^\theta)f.
$$
Thus, all that is left to check, is that integrating the left-hand side yields $\Delta_G^{-1/4}\star_\theta (\lambda(f)\eta)\star_\theta\Delta_G^{1/4}$, that is, that conjugation by $\Delta_G^{-1/4}$ with respect to the product $\star_\theta$ commutes with integration. But this is clear, as conjugation by $\Delta_G^{-1/4}$ is a continuous map on $\Sch(G)$.
\ep

\bigskip

\section{Deformations of \texorpdfstring{C$^*$}{C*}-algebras}\label{sec:s3}

We now generalize the results of the previous section to actions of K\"ahlerian Lie groups with negative curvature on C$^*$-algebras.

\subsection{Elementary case}
We start with the case of an elementary K\"ahlerian Lie group $G$. Consider a $C^*$-algebra $A$. Then one can define
in a straightforward way the $A$-valued versions   of $\Sch(G)$ and~$\B(G)$. Since $\Sch(G)$
is nuclear as a locally convex topological vector space, one can also define $\Sch(G,A)$ as the unique
completion of the algebraic tensor product $\Sch(G)\otimes_{\rm alg}A$.

As shown in \cite{BG}, the oscillatory integrals and the product $\star_\theta$ make sense for $A$-valued functions, so $\star_\theta$ is defined on $\B(G,A)$ by the same formula~\eqref{SP}. Furthermore, $(\Sch(G,A),\star_\theta)$ is an ideal in $(\B(G,A),\star_\theta)$ and
we have a representation $\pi_\theta\otimes\iota$ of $(\B(G,A),\star_\theta)$ on the Hilbert $A$-module $L^2(G)\otimes A$ defined by
$$
(\pi_\theta\otimes\iota)(f)\xi=f\star_\theta\xi\ \ \text{for}\ \ \xi\in\Sch(G,A)\subset L^2(G)\otimes A.
$$
By \cite[Theorems 8.20 \& 8.33]{BG} this representation satisfies the same estimate \eqref{to-get-sunlight}.

It is almost a tautological statement that the representation $\pi_\theta\otimes\iota$ of $(\B(G,A),\star_\theta)$ can be described in terms of $\pi_\theta$ and the slice (or Fubini) maps $\iota\otimes\nu\colon M(\K\otimes A)\to M(\K)=\B(L^2(G))$ for $\nu\in A^*$. Namely, we have the following result.

\begin{lemma}\label{l:os-free2}
For any $f\in\B(G,A)$, the operator $(\pi_\theta\otimes\iota)(f)\in M(\K\otimes A)$ is the unique element satisfying
$$
(\iota\otimes\nu)\big((\pi_\theta\otimes\iota)(f)\big)=\pi_\theta\big((\iota\otimes\nu)(f)\big)\ \ \text{for all}\ \nu\in A^*.
$$
\end{lemma}

\bp
Uniqueness is clear. Replacing, if necessary, $A$ by its unitization, we may assume that $A$ is unital.
Take $\xi\in\Sch(G)$. Then
$$
(\pi_\theta\otimes\iota)(f)(\xi\otimes 1)=f\star_\theta(\xi\otimes 1).
$$
By definition, the oscillatory integrals commute with the slice maps. Hence, applying $\iota\otimes\nu$ to the above equality, we get
$$
(\iota\otimes\nu)\big((\pi_\theta\otimes\iota)(f)\big)\xi=(\iota\otimes\nu)(f)\star_\theta\xi,
$$
which is what we need.
\ep

\begin{remark}
The above lemma together with Lemma~\ref{l:os-free} show that $\pi_\theta\otimes\iota$ can be described entirely in terms of ordinary integrals and slice maps.
\end{remark}

Consider now a strongly continuous action $\alpha$ of $G$ on a C$^*$-algebra $A$. For $a\in A$ we denote by
$\check\alpha(a)$ the $A$-valued right uniformly continuous function $g\mapsto\alpha_g(a)$. Denote also by $A^\infty$ the Fr\'echet algebra of smooth elements in $A$, which is exactly the set of all elements $a\in A$ such that $\check\alpha(a)\in \B(G,A)$. Then a new product $\star_\theta$ can be defined on $A^\infty$ by
$$
\check\alpha(a\star_\theta b)=\check\alpha(a)\star_\theta\check\alpha(b).
$$
The deformation $A_\theta$ of  $A$ is defined in \cite[Section~8.5]{BG} by
$$
A_\theta=[(\pi_\theta\otimes\iota)\check\alpha(a): a\in A^\infty]\subset M(\K\otimes A).
$$
Once again we remark that this definition of~$A_\theta$ is slightly different from the one given in \cite{BG},
but equivalent to it, as we use the representation~$\pi_\theta$ on $L^2(G)$ instead of a quasi-equivalent irreducible
representation on $L^2(Q)$.

On the other hand, a deformation of $A$ can be defined using the quantization maps studied in Section~\ref{s:qm}. Namely, following~\cite[Section~3.2]{NTdef2}, let
$$
A_{\Omega_\theta}=[(T_\nu\otimes\iota)\alpha(a): a\in A,\ \nu\in \K^*]\subset M(\K\otimes A),
$$
where $\alpha(a)$ is the $A$-valued function $g\mapsto\alpha_{g^{-1}}(a)$. In general it is apparently necessary to consider the algebra generated by the elements $(T_\nu\otimes\iota)\alpha(a)$ before taking the norm closure, but since $\Omega_\theta$ is regular by Corollary~\ref{cregular}, we don't have to do this in the present case by \cite[Theorem~3.7]{NTdef2}. This also follows from the proof of the next theorem.

\begin{theorem}
\label{main-stuff}
For any C$^*$-algebra $A$ equipped with a strongly continuous action $\alpha$ of an elementary
K\"ahlerian Lie group $G$, we have
$A_{\Omega_\theta}=(\bj\otimes 1)A_\theta(\bj\otimes1)\subset M(\K\otimes A)$.
\end{theorem}

\bp
Define maps $\check T_\nu$ by $\check T_\nu(f)=T_\nu(\check f)$ for $f\in L^\infty(G)$. We can then write
$$
A_{\Omega_\theta}=[(\check T_{\omega_{f_1,f_2}}\otimes\iota)\check\alpha(a): a\in A,\ f_1,f_2\in
\Sch(G)].
$$

We claim that
$$
\Ad(\bj\otimes1)(\check T_{\omega_{f_1,f_2}}\otimes\iota)=(\pi_\theta\otimes\iota)
\mathcal R(\eta^\theta_{f_1,f_2})
\ \ \text{on}\ \ C_b(G,A),
$$
where $\mathcal R_g$ and $\RR(\eta)$ are defined in the same way as before, but now on $A$-valued functions, so $(\RR_gf)(h)=f(hg)$ and $\RR(\eta)=\int\eta(g)\,\RR_g\,dg$. First of all, note that the same argument as in the proof of Lemma~\ref{something-in-my-eye} shows that $\RR(\eta^\theta_{f_1,f_2})$ maps $C_b(G,A)$ into $\B(G,A)$, so both sides of the above identity are at least well-defined. In order to prove the identity it suffices to show that we get the same operators if we apply $\iota\otimes\nu$ to both sides for all $\nu\in A^*$. When we apply $\iota\otimes\nu$, the left-hand side gives $(\Ad \bj)\check T_{\omega_{f_1,f_2}}(\iota\otimes\nu)$. On the right-hand side, using Lemma~\ref{l:os-free2} and that $\RR(\eta)$ commutes with the slice maps, we get $\pi_\theta \mathcal R(\eta^\theta_{f_1,f_2})(\iota\otimes\nu)$.
Since $(\Ad \bj)\check T_{\omega_{f_1,f_2}}=\pi_\theta \mathcal R(\eta^\theta_{f_1,f_2})$ on $C_b(G)$ by Proposition~\ref{pquantmaps}, our claim is therefore proved.

Since for $a\in A$ we have $\check\alpha(a)\in C_b(G,A)$, it follows that
$$
(\bj\otimes 1)A_{\Omega_\theta}(\bj\otimes 1)=[(\pi_\theta\otimes\iota)\mathcal R(\eta^\theta_{f_1,f_2})\check\alpha(a): a\in A,\ f_1,f_2\in \Sch(G)].
$$
Therefore it remains to show that the right-hand side of this identity coincides with $A_\theta$.
In order to see this, for $\eta\in L^1(G)$ consider the operator $\RR^\alpha(\eta)$ on $A$ defined by
$$
\RR^\alpha(\eta)a=\int_G\eta(g)\,\alpha_g(a)\,dg.
$$
Then $\mathcal R(\eta)\check\alpha(a)=\check\alpha(\RR^\alpha(\eta)a)$ and therefore we must check that
$$
[(\pi_\theta\otimes\iota)\check\alpha(\RR^\alpha(\eta^\theta_{f_1,f_2})a): a\in A\, f_1,f_2\in \Sch(G)]=A_\theta.
$$
For this it suffices to show that the elements $\RR^\alpha(\eta^\theta_{f_1,f_2})a$ for $a\in A$ and
$f_1,f_2\in\Sch(G)$ span a dense subspace of the Fr\'echet space $A^\infty$.

By definition, we have
$$
\RR^\alpha(\eta^\theta_{f_1,f_2})a=\int_G \Delta_G^{-1}(g)\, (\bar f_2\star_\theta f_1)(g^{-1})\,
\alpha_g(a)\,dg.
$$
By \cite[Proposition 5.19]{BG} we have a bounded approximate unit
for the Fr\'echet algebra $(\Sch(G),\star_\theta)$. Letting $\bar f_2$ run through such
an approximate unit shows that the closure (with respect to the topology on $A^\infty$) of linear
combinations of all elements of $A^\infty$ of the form $\RR^\alpha(\eta^\theta_{f_1,f_2})a$ contains
all linear combinations of elements of the form $\RR^\alpha(f)a$ with $f\in \Sch(G)$ and $a\in A^\infty$. Hence
it contains the G{\aa}rding subspace of $A$, which by the Dixmier-Malliavin theorem coincides with~ $A^\infty$.
\ep

\subsection{General case}

Consider now an arbitrary negatively curved K\"ahlerian Lie group $G$ with Pyatetskii-Shapiro
decomposition~\eqref{KLG}.
For $g=g_1\dots g_n$ and $g'=g'_1\dots g'_n\in G$ with $g_i,g_i'\in G_i$, define a two-point kernel on $G$ by
\begin{equation}
\label{very-cute}
K_\theta(g,g'):=\prod_{i=1}^n K^{G_i}_\theta(g_i,g'_i),
\end{equation}
where $K^{G_i}_\theta(g_i,g'_i)$ is the two-point kernel on the elementary K\"ahlerian Lie group
$G_i$ as defined in~\eqref{KTheta}. Let also $A$ be a $C^*$-algebra endowed with a strongly continuous
action $\alpha$ of $G$. Then, as shown in~\cite {BG}, one can define a deformation  $A_\theta$ of $A$
exactly as we did in the case of an elementary K\"ahlerian Lie group. We also have a dual unitary $2$-cocycle on $G$ defined by
$$
\Omega_\theta=\int_{G\times G}\overline{K_\theta(g_1,g_2)}\, \lambda_{g_1^{-1}}\otimes\lambda_{g_2^{-1}}
\,dg_1\,dg_2,
$$
and hence can define a deformation $A_{\Omega_\theta}$ of $A$.

\begin{theorem}
\label{main-stuff2}
For any C$^*$-algebra $A$ equipped with a strongly continuous action $\alpha$ of a negatively curved
K\"ahlerian Lie group $G$, we have
$A_{\Omega_\theta}=(\bj\otimes 1)A_\theta(\bj\otimes1)\subset M(\K\otimes A)$.
\end{theorem}

\bp Since the action  of $\big(\dots\big(G_n\ltimes G_{n-1}\big)\ltimes\dots\big)\ltimes G_{i+1}$ on $G_i$ leaves the Haar measure and the kernel~$K^{G_i}_\theta$ invariant, the modular function of $G$ is the product  of the modular functions of the factors, and the kernel $K_\theta$ still satisfies identity \eqref{weird-inv}.
Routine verifications show then that all the previous arguments extend with only minor changes to general negatively curved
K\"ahlerian Lie groups. In fact, the only places, where we used that we worked with elementary K\"ahlerian Lie groups, were Lemmas~\ref{l:key} \& \ref{trace-of-foxes} and Proposition~\ref{pmodular}. They can be easily extended to the general case by working with the dense subalgebra $\Sch(G_1)\otimes_{\mathrm{alg}}\dots \otimes_{\mathrm{alg}}\Sch(G_n)$ of $\Sch(G)$. We leave the details to the reader.
\ep

\bigskip

\appendix

\section{Decomposition of the dual cocycle}\label{ap:A}

In this appendix we obtain an explicit expression for $\Omega_\theta$ for an elementary K\"ahlerian Lie group $G=Q\ltimes P$, showing that this dual cocycle is well-defined and coisometric. (The main part of the paper wrongly states that the cocycle is unitary, see Appendix~\ref{ap:B} for further discussion.)

\smallskip

Define an operator $\F_{P,\theta}\colon L^2(G)\to L^2(G)$ by
$$
(\F_{P,\theta} f)(q,p)=(\pi\theta)^{-(d+1)/2}\int_{P} e^{\tfrac{2i}\theta p.p'}f(q,p')\,dp',
$$
so, up to rescaling $p$ by $\theta$ and normalization, $\F_{P,\theta}$ is the partial Fourier transform associated to the maximal Abelian subgroup $P\simeq\R^{d+1}$ of $G$. It is clearly unitary.

Recall that we denote the coordinates of $G$ by $(a,v,t)$. We will write $v\in\R^{2d}$ as $v=(n,m)$ with $n,m\in\R^d$, so that $\omega_0(v,v')=n.m'-m.n'$. Define a map  $\Phi\colon G\times G\to G\times G$ by
\begin{align*}
\Phi(g_1,g_2):=
\big(&a_1-\tfrac12\arcsinh(e^{-2a_2}t_2),\tfrac{\alpha^2}{\alpha'}n_1+\tfrac{t_2}{2\alpha'}
m_1+\tfrac{\alpha t_2}{2\alpha'} n_2-\tfrac\alpha{\alpha'}m_2,\\
&\qquad\qquad\qquad\qquad\qquad\ \tfrac{t_1^2t_2}{8\alpha'}n_1+\tfrac{\alpha^2}{\alpha'}
m_1-\tfrac{\alpha t_1 t_2}{4\alpha'} n_2+\tfrac{\alpha t_1}{2\alpha'}m_2,t_1;\\
&a_2+\tfrac12\arcsinh(e^{-2a_1}t_1),-\tfrac{\alpha t_1}{2\alpha'}n_1+
\tfrac\alpha{\alpha'}m_1+\tfrac{\alpha^2}{\alpha'}n_2+\tfrac{t_1}{2\alpha'}m_2,\\
&\qquad\qquad\qquad\qquad\qquad\ \tfrac{\alpha t_1 t_2}{4\alpha'}n_1-\tfrac{\alpha t_2}{2\alpha'}
m_1+\tfrac{t_1 t_2^2}{8\alpha'} n_2+\tfrac{\alpha^2}{\alpha'}m_2,t_2\big),
\end{align*}
where
\begin{equation}\label{eq:alpha-numbers}
\alpha=\tfrac12 \big(e^{2a_1}+(e^{4a_1}+t_1^2)^{1/2}\big)^{1/2}\big(e^{2a_2}+(e^{4a_2}+t_2^2)^{1/2}\big)^{1/2},\ \
\alpha'=\alpha^2+\tfrac{1}{4}t_1t_2.
\end{equation}
Consider the corresponding operator $U_\Phi$ on $L^2(G\times G)$:
$$
(U_\Phi\varphi)(g_1,g_2):=\Jac_\Phi(g_1,g_2)^{1/2}\varphi(\Phi(g_1,g_2)).
$$

\begin{theorem}\label{thm:A}
For any $\theta\in\R^*$, we have
$$
\Omega_\theta=(\F_{P,\theta}^{-1}\otimes \F_{P,\theta}^{-1})U_\Phi(\F_{P,\theta}\otimes\F_{P,\theta})\ \ \text{on}\ \ L^2(G\times G).
$$
\end{theorem}

\bp By definition, $\Omega_\theta$ acts on a vector $\vf\in\Sch(G\times G)$ as
\begin{align*}
(\Omega_\theta \vf)(g_1,g_2)&=\int_{G\times G} \overline{K_\theta(g_1',g_2')} \, \vf(g_1'g_1,g_2'g_2)\,
dg_1'\,dg_2'\\
&=\Delta_G^{-1}(g_1)\Delta_G^{-1}(g_2)\int_{G\times G} \overline{K_\theta(g_1'g_1^{-1},g_2'g_2^{-1})} \, \vf(g_1',g_2')\,dg_1'\,dg_2'.
\end{align*}
Hence, the distributional kernel of the operator $\Omega_\theta$ is given by
\begin{align*}
[\Omega_\theta](g_1,g_2;g_1',g_2')&=\Delta_G^{-1}(g_1)\Delta_G^{-1}(g_2)
\overline{K_\theta(g_1'g_1^{-1},g_2'g_2^{-1})},
\end{align*}
which in coordinates reads as
\begin{align*}
&\frac{4}{(\pi\theta)^{2d+2}}
e^{(2d+2)(a_1+a_2)}
\big(\cosh(a_1'-a_1)\cosh(a_2'-a_2)\cosh(a_1'-a_1-a_2'+a_2)\big)^d\\
&\times\big(\cosh(2a_1'-2a_1)\cosh(2a_2'
-2a_2)\cosh(2a_1'-2a_1-2a_2'+2a_2)\big)^{1/2}\\
&\times\exp\big\{-\tfrac{2i}\theta\sinh(2a_1'-2a_1)e^{2a_2}(t_2'-t_2-\tfrac12(n_2'.m_2-n_2.m_2'))\big\}\\
&\times\exp\big\{\tfrac{2i}\theta\sinh(2a_2'-2a_2)e^{2a_1}(t_1'-t_1-\tfrac12(n_1'.m_1-n_1.m_1'))\big\}\\
&\times\exp\big\{-\tfrac{2i}\theta\cosh(a_1'-a_1)\cosh(a_2'-a_2)e^{a_1+a_2}\big((n_1'-n_1).(m_2'-m_2)-
(n_2'-n_2).(m_1'-m_1)\big)\big\}.
\end{align*}
On the other hand, the distributional kernel of $\F_{P,\theta}$ is given by
$$
[\F_{P,\theta}](g,g')=(\pi\theta)^{-(d+1)/2}\delta_{(a,n)}(a',n')\exp\{\tfrac{2i}\theta(tt'+m.m')\}.
$$
From this we deduce the following expression for the kernel of $(\F_{P,\theta}\otimes \F_{P,\theta})\Omega_\theta
(\F_{P,\theta}^{-1}\otimes \F_{P,\theta}^{-1})$:
\begin{align*}
&[(\F_{P,\theta}\otimes \F_{P,\theta})\Omega_\theta
(\F_{P,\theta}^{-1}\otimes \F_{P,\theta}^{-1})](g_1,g_2;g_3,g_4)\\
&=\int [\F_{P,\theta}\otimes \F_{P,\theta}](g_1,g_2;g_1',g_2')
[\Omega_\theta](g_1',g_2';g_1'',g_2'')[\F_{P,\theta}^{-1}\otimes \F_{P,\theta}^{-1}](g_1'',g_2'';g_3,g_4)\,dg_1'\,dg_2'\,dg_1''\,dg_2''
\\
&=\frac{4}{(\pi\theta)^{4d+4}}e^{(2d+2)(a_1+a_2)}
\big(\cosh(a_3-a_1)\cosh(a_4-a_2)\cosh(a_3-a_1-a_4+a_2)\big)^d\\
&\quad\times\big(\cosh(2a_3-2a_1)\cosh(2a_4
-2a_2)\cosh(2a_3-2a_1-2a_4+2a_2)\big)^{1/2}\\
&\quad\times\int \exp\{\tfrac{2i}\theta t_1'(t_1-\sinh(2a_4-2a_2)e^{2a_1})\}
\exp\{\tfrac{2i}\theta t_2'(t_2+\sinh(2a_3-2a_1)e^{2a_2})\}\\
&\quad\times\exp\{-\tfrac{2i}\theta t_1''(t_3-\sinh(2a_4-2a_2)e^{2a_1})\}
\exp\{-\tfrac{2i}\theta t_2''(t_4+\sinh(2a_3-2a_1)e^{2a_2})\}\\
&\quad\times\exp\{\tfrac{2i}\theta m_1'.(m_1-\tfrac12\sinh(2a_4-2a_2)e^{2a_1}n_3-\cosh(a_3-
a_1)\cosh(a_4-a_2)e^{a_1+a_2}(n_4-n_2))\}\\
&\quad\times\exp\{\tfrac{2i}\theta m_2'.(m_2+\tfrac12\sinh(2a_3-2a_1)e^{2a_2}n_4 +\cosh(a_3-
a_1)\cosh(a_4-a_2)e^{a_1+a_2}(n_3-n_1))\}\\
&\quad\times\exp\{-\tfrac{2i}\theta m_1''.(m_3-\tfrac12\sinh(2a_4-2a_2)e^{2a_1}n_1-
\cosh(a_3-a_1)\cosh(a_4-a_2)e^{a_1+a_2}(n_4-n_2))\} \\
&\quad\times\exp\{-\tfrac{2i}\theta m_2''.(m_4+\tfrac12\sinh(2a_3-2a_1)e^{2a_2}n_2+
\cosh(a_3-a_1)\cosh(a_4-a_2)e^{a_1+a_2}(n_3-n_1))\}\\
&\quad\times dt_1'\,dt_2'\,dt_1''\,dt_2''\,dm_1'\,dm_2'\,dm_1''\,dm_2''.
\end{align*}
Integrating out the phase factors produces delta-factors and we get that the above expression equals
\begin{align}
&4e^{(2d+2)(a_1+a_2)}
\big(\cosh(a_3-a_1)\cosh(a_4-a_2)\cosh(a_3-a_1-a_4+a_2)\big)^d\nonumber\\
&\times\big(\cosh(2a_3-2a_1)\cosh(2a_4
-2a_2)\cosh(2a_3-2a_1-2a_4+2a_2)\big)^{1/2}\delta_0(\Xi(g_1,g_2,g_3,g_4)),\label{eq:kernel1}
\end{align}
where $\Xi\colon G\times G\times G\times G\to\R^{4d+4}$ is the map defined by
\begin{multline*}
\Xi(g_1,g_2,g_3,g_4):=(t_3-t_1,t_4-t_2,\\
t_1-\sinh(2a_4-2a_2)e^{2a_1},t_2+\sinh(2a_3-2a_1)e^{2a_2},
 A\begin{pmatrix}n_1 \\ m_1 \\ n_2 \\ m_2\end{pmatrix}
+B\begin{pmatrix}n_3 \\ m_3 \\ n_4 \\ m_4\end{pmatrix}),
\end{multline*}
with
$\displaystyle
A=\begin{pmatrix} 0 & 1 & a & 0\\-a & 0 & 0 & 1\\-b & 0 & a & 0\\ -a & 0 & c & 0\end{pmatrix}$,
$\displaystyle B=\begin{pmatrix}-b & 0 & -a & 0\\a & 0 & c & 0\\0 & 1 & -a & 0\\ a & 0 & 0 & 1\end{pmatrix}$ and
$$
a=\cosh(a_3-a_1)\cosh(a_4-a_2)e^{a_1+a_2},\ \
b=\tfrac12\sinh(2a_4-2a_2)e^{2a_1}\ \
c=\tfrac12\sinh(2a_3-2a_1)e^{2a_2}.
$$
It is not difficult to check that
\begin{equation}\label{eq:det}
\det(A)=\det(B)=bc-a^2=-e^{2(a_1+a_2)}\cosh(a_3-a_1)\cosh(a_4-a_2)\cosh(a_3-a_1-a_4+a_2)
\end{equation}
and that for $a_3=a_1-\tfrac12\arcsinh(e^{-2a_2}t_2)$ and $a_4=a_2+\tfrac12\arcsinh(e^{-2a_1}t_1)$ we have
$$
a=\alpha,\ \ b=\frac{t_1}{2},\ \ c=-\frac{t_2}{2},\ \ a^2-bc=\alpha',
$$
where  $\alpha$ and $\alpha'$ are defined by~\eqref{eq:alpha-numbers}. From this we easily conclude that for fixed $(g_1,g_2)$ the only solution of the equation $\Xi(g_1,g_2,g_3,g_4)=0$ is $(g_3,g_4)=\Phi(g_1,g_2)$. Observe also that $\Xi(g_1,g_2,\cdot,\cdot)$ is a diffeomorphism of $\R^{4d+4}$ onto itself. It follows that
$$
\delta_0(\Xi(g_1,g_2,g_3,g_4))=\Jac_{\Xi(g_1,g_2,\cdot,\cdot)}(g_3,g_4)^{-1}\delta_{\Phi(g_1,g_2)}(g_3,g_4).
$$
In view of~\eqref{eq:det} we have
\begin{multline}
\Jac_{\Xi(g_1,g_2,\cdot,\cdot)}(g_3,g_4)=4e^{(2d+2)(a_1+a_2)}\cosh(2a_3-2a_1) \cosh(2a_4-2a_2)\\
\times\big(\cosh(a_3-a_1)\cosh(a_4-a_2)\cosh(a_3-a_1-a_4+a_2)\big)^d.\label{eq:JacXi}
\end{multline}
Therefore~\eqref{eq:kernel1} equals
$$
\frac{\cosh(2a_3-2a_1-2a_4+2a_2)^{1/2}}{\big(\cosh(2a_3-2a_1)\cosh(2a_4
-2a_2)\big)^{1/2}}\delta_{\Phi(g_1,g_2)}(g_3,g_4).
$$

To finish the proof of the equality $(\F_{P,\theta}\otimes \F_{P,\theta})\Omega_\theta
(\F_{P,\theta}^{-1}\otimes \F_{P,\theta}^{-1})=U_\Phi$ it remains to check that
$$
\frac{\cosh(2a_3-2a_1-2a_4+2a_2)}{\cosh(2a_3-2a_1)\cosh(2a_4-2a_2)}\Big|_{(g_3,g_4)=\Phi(g_1,g_2)}=\Jac_\Phi(g_1,g_2).
$$
Since $\Xi(g_1,g_2,\Phi(g_1,g_2))=0$, we have
$$
\Jac_\Phi(g_1,g_2)=\frac{\Jac_{\Xi(\cdot,\cdot,\Phi(g_1,g_2))}(g_1,g_2)}{\Jac_{\Xi(g_1,g_2,\cdot,\cdot)}(\Phi(g_1,g_2))}.
$$
A straightforward computation using~\eqref{eq:det} yields
\begin{multline*}
\Jac_{\Xi(\cdot,\cdot,g_3,g_4)}(g_1,g_2)=4e^{(2d+2)(a_1+a_2)}\cosh(2a_3-2a_1-2a_4+2a_2)\\
\times\big(\cosh(a_3-a_1)\cosh(a_4-a_2)\cosh(a_3-a_1-a_4+a_2)\big)^d.
\end{multline*}
Together with~\eqref{eq:JacXi} this gives the required identity.
\ep

\bigskip

\section{Erratum}\label{ap:B}

The paper relies in a crucial way on unitarity of the cocycle $\Omega_\theta$, but it turns out that $\Omega_\theta$ is only a coisometry. This can be seen from Theorem~\ref{thm:A}: the map~$\Phi$ is injective but not surjective (see Lemma~\ref{lem:O} below). In fact, it seems doubtful that already the simplest K\"ahlerian Lie group $G$, the connected component of the $ax+b$ group over the reals, supports any nontrivial dual unitary $2$-cocycle. At the very least such a cocycle must be such that the corresponding deformation of the function algebra $L^\infty(G)$ is not a type I factor, see~\cite[Remark~2.12]{BGNT2}. The entire nonconnected $ax+b$ group does support a nontrivial dual cocycle though~\cite{BGNT2}.

The unitarity of $\Omega_\theta$ was stated for the first time in~\cite{NTdef2} as a consequence of the identity from~\cite[Proposition~8.47]{BG}, which therefore is also wrong. The origin of the mistake is similar: a change of variables closely related to $\Phi$ is not bijective, see erratum to~\cite{BG}.

Although the dual cocycle $\Omega_\theta$ is not unitary, it still has some nice properties that allow one to develop a deformation theory similar to~\cite{NTdef2}. Within the framework of this extended theory the main results of the paper remain true as stated, with essentially identical proofs. This is what we are going to explain in this erratum.

\subsection{Deformation by nonunitary cocycles}\label{sec:def}

We follow the conventions of~\cite{NTdef2}. Let $G$ be a locally compact quantum group and $\hat W$ be the multiplicative unitary of the dual quantum group. Assume we are given an element $\Omega\in L^\infty(\hat G)\bar\otimes L^\infty(\hat G)$ satisfying the cocycle identity
$$
(\Omega\otimes1)(\Dhat\otimes\iota)(\Omega)=(1\otimes\Omega)(\iota\otimes\Dhat)(\Omega).
$$
We assume in addition that:
\begin{enumerate}
\item[(1)] the element $\Omega$ is coisometric: $\Omega\Omega^*=1$;
\item[(2)] there exists a unitary $X\in L^\infty(\hat G)$ such that $(XJ)^2=1$ and $\Omega=(X\otimes X)(\hat R\otimes\hat R)(\Omega^*_{21})\Dhat(X)^*$;
\item[(3)] we have $(\Dhat\otimes\iota)(\Omega)(\iota\otimes\Dhat)(\Omega^*)=(\Omega^*\otimes1)(1\otimes\Omega)$.
\end{enumerate}
Recall that $\hat R$ denotes the unitary antipode on $L^\infty(\hat G)$, so that $\hat R(x)=Jx^*J$.

We then let
$$
\tilde J:=XJ,\qquad L:=(\tilde J\otimes \hat J)\Omega\hat W^*(J\otimes\hat J)
$$
and define, for $\nu\in\K^*=\B(L^2(G))_*$,  a \emph{quantization map}
$$
T_\nu\colon L^\infty(G)\to\B(L^2(G))\ \ \text{by}\ \ T_\nu(x):=(\iota\otimes\nu)(L(x\otimes1)L^*).
$$
Given a left action $\alpha\colon A\to M(C_0(G)\otimes A)$ on a C$^*$-algebra $A$, we define the \emph{$\Omega$-deformation} of~$A$ as the C$^*$-subalgebra $A_\Omega\subset M(\K\otimes A)$ generated by the elements $(T_\nu\otimes\iota)\alpha(a)$ for $a\in A$ and $\nu\in\K^*$.
Note that although this definition depends on the unitary $X$, a different choice would give an isomorphic C$^*$-algebra.

The definition of $A_\Omega$ as such does not require conditions (1)--(3): we could take any unitary in $L^\infty(\hat G)$ for $X$ in order to define $L$.
However, without additional assumptions it is not clear how reasonable this definition is, since already for the function algebra there is another more transparent definition. Namely, we can define a product $\star_\Omega$ on the Fourier algebra $A(G):=L^\infty(\hat G)_*$ of~$G$~by
$$
\omega\star_\Omega\nu=(\omega\otimes\nu)(\Dhat(\cdot)\Omega^*).
$$
It is straightforward to check that the cocycle identity for $\Omega$ can be written as
$$
(\Dhat\otimes\iota)(\hat W\Omega^*)\Omega^*_{12}=(\hat W\Omega^*)_{13}(\hat W\Omega^*)_{23},
$$
see~\cite[Section~2.1]{NTdef2}. This implies that we have a representation
$$
\pi_\Omega\colon (A(G),\star_\Omega)\to \B(L^2(G))\ \ \text{defined by}\ \ \pi_\Omega(\omega):=(\omega\otimes\iota)(\hat W\Omega^*).
$$
How is the algebra $\pi_\Omega(A(G))$ related to the deformation $C_0(G)_\Omega$ of $C_0(G)$ with respect to the left action of $G$ on itself by translations?
This can be answered under conditions (1)--(3) thanks to the following result, cf.~identity~(3.1) in~\cite{NTdef2}.

\begin{proposition}\label{prop:eq3.1}
If an element $\Omega\in L^\infty(\hat G)\bar\otimes L^\infty(\hat G)$ satisfies conditions {\rm (2)} and {\rm (3)}, then
\begin{equation}\label{eq:3.1}
L_{23}\hat W_{12}L^*_{23}=(\hat W\Omega^*)_{12}L_{13}.
\end{equation}
\end{proposition}

\bp
By the definition of $L$ we have to prove that
\begin{multline*}
(J\otimes\tilde J\otimes \hat J)\Omega_{23}\hat W^*_{23}(J\otimes J\otimes\hat J)\hat W_{12}
(J\otimes J\otimes \hat J)\hat W_{23}\Omega_{23}^*(J\otimes \tilde J\otimes\hat J)\\
=\hat W_{12}\Omega^*_{12}(\tilde J\otimes J\otimes \hat J)\Omega_{13}\hat W^*_{13}(J\otimes J\otimes\hat J).
\end{multline*}
As $(J\otimes\hat J)\hat W=\hat W^*(J\otimes\hat J)$ and $\hat W_{23}\hat W_{12}\hat W^*_{23}=\hat W_{12}\hat W_{13}$, this is equivalent to
$$
(J\otimes\tilde J\otimes \hat J)\Omega_{23}(J\otimes J\otimes\hat J)\hat W_{12}\hat W_{13}(J\otimes J\otimes\hat J)\Omega_{23}^*(1\otimes X\otimes1)
=\hat W_{12}\Omega^*_{12}(\tilde J\otimes J\otimes \hat J)\Omega_{13}\hat W^*_{13},
$$
hence to
$$
(J\otimes\tilde J\otimes \hat J)\Omega_{23}(J\otimes J\otimes\hat J)\hat W_{12}(J\otimes J\otimes\hat J)\hat W_{13}^*\Omega_{23}^*(1\otimes X\otimes1)
=\hat W_{12}\Omega^*_{12}(\tilde J\otimes J\otimes \hat J)\Omega_{13}\hat W^*_{13}.
$$
Multiplying by $\hat W^*_{12}(1\otimes1\otimes J\hat J)$ on the left and by $\hat W_{13}(J\otimes X^*J\otimes J)$ on the right and using that $XJ=JX^*$, we get
\begin{multline*}
\hat W_{12}^*(J\otimes\tilde J\otimes J)\Omega_{23}(J\otimes J\otimes\hat J)\hat W_{12}(J\otimes J\otimes\hat J)\hat W_{13}^*\Omega_{23}^*\hat W_{13}(J\otimes J\otimes J)\\
=\Omega^*_{12}(\tilde J\otimes XJ\otimes J)\Omega_{13}(J\otimes J\otimes J).
\end{multline*}
Since we can replace $\hat J$ by $J$ in this expression, we get
$$
\hat W_{12}^*(1\otimes X \otimes 1)(\hat R\otimes\hat R)(\Omega^*)_{23}\hat W_{12}(\hat R\otimes \hat R\otimes \hat R)(\hat W_{13}^*\Omega_{23}\hat W_{13})
=\Omega^*_{12}(X\otimes X\otimes1)(\hat R\otimes\hat R)(\Omega^*)_{13}.
$$
Recalling that the coproduct on $L^\infty(\hat G)$ is given by $\Dhat(x)=\hat W^*(1\otimes x)\hat W$, we obtain
$$
\Dhat(X)_{12}(\Dhat\otimes\iota)(\hat R\otimes\hat R)(\Omega^*)(\hat R\otimes\hat R\otimes\hat R)(\iota\otimes\Dhat)(\Omega)_{213}
=\Omega^*_{12}(X\otimes X\otimes1)(\hat R\otimes\hat R)(\Omega^*)_{13}.
$$
Applying $\hat R\otimes\hat R\otimes\hat R$ and flipping the first two factors, we then get
$$
(\iota\otimes\Dhat)(\Omega)(\Dhat\otimes\iota)(\Omega^*)\Dhat(X)_{12}
=\Omega^*_{23}(X\otimes X\otimes1)(\hat R\otimes\hat R)(\Omega^*)_{21},
$$
where we used that $(\hat R\otimes\hat R)\Dhat=\Dhat^{\mathrm{op}}\hat R$ and $\hat R(X)=X$. By virtue of assumption (2) this is equivalent~to
$$
(\iota\otimes\Dhat)(\Omega)(\Dhat\otimes\iota)(\Omega^*)
=\Omega^*_{23}\Omega_{12},
$$
which is condition (3).
\ep

\begin{corollary}[cf.~{\cite[Proposition~3.1]{NTdef2}}]
If $\Omega$ is a dual cocycle satisfying conditions {\rm (1)--(3)}, then
\begin{equation}\label{eq:function-alg}
\overline{\pi_\Omega(A(G))}=[T_\nu(x)\colon x\in C_0(G),\ \nu\in \K^*].
\end{equation}
It follows that $\overline{\pi_\Omega(A(G))}$ is a nondegenerate C$^*$-algebra of operators on $L^2(G)$ and
$$
C_0(G)_\Omega=V(\overline{\pi_\Omega(A(G))}\otimes1) V^*,
$$
where $V=(\hat J\otimes\hat J)\hat W(\hat J\otimes\hat J)$.
\end{corollary}

We remind that the square brackets denote the norm closure of the linear span.

\bp
Slicing the first and the third legs of~\eqref{eq:3.1} and remembering that $L_{13}$ is coisometric we immediately get the first statement. Next, the right hand side of~\eqref{eq:function-alg} is a self-adjoint space of operators, since the maps $T_\nu$ are completely positive for positive $\nu$. On the other hand, the left hand side is an algebra. Hence both sides give a C$^*$-algebra. This C$^*$-algebra acts nondegenerately on~$L^2(G)$, because $C_0(G)$ acts nondegenerately and $T_\nu(1)=\nu(1)1$.

Since the first leg of $L$ is in $L^\infty(\hat G)$ while the first leg of $V$ is in $L^\infty(\hat G)'$, we have
$$
V(T_\nu(x)\otimes 1)V^*=(T_\nu\otimes\iota)(V(x\otimes1)V^*)=(T_\nu\otimes\iota)\Delta(x)
$$
for $x\in C_0(G)$. This proves the last statement of the corollary.
\ep

\subsection{Dual cocycles on K\"{a}hlerian Lie groups}

As in~Section~\ref{sec:s3}, let $G$ be a K\"{a}hlerian Lie group, $\theta\in\R^*$ and $\Omega_\theta$ be the dual cocycle defined by the kernel $K_\theta$. We want to check that $\Omega_\theta$ satisfies conditions (1)--(3). The operator $\Omega_\theta$ is coisometric by Theorem~\ref{thm:A} and injectivity of the map $\Phi$ there (see Lemma~\ref{lem:O} below). Condition (2) is also satisfied for $X=1$, since $\overline{K_\theta(g,h)}=K_\theta(h,g)$. As for (3), the cocycle identity implies that
$$
(\Omega_\theta^*\Omega_\theta\otimes1)(\Dhat\otimes\iota)(\Omega_\theta)(\iota\otimes\Dhat)(\Omega_\theta^*)=(\Omega_\theta^*\otimes1)(1\otimes\Omega_\theta).
$$
It follows that if we let $F_\theta:=\Omega_\theta^*\Omega_\theta$ then condition (3) is equivalent to
$$
(F_\theta\otimes1)(\Dhat\otimes\iota)(\Omega_\theta)(\iota\otimes\Dhat)(F_\theta)=(\Dhat\otimes\iota)(\Omega_\theta)(\iota\otimes\Dhat)(F_\theta).
$$
We will show that the following stronger property holds.

\begin{proposition}\label{prop:mystery}
We have $(F_\theta\otimes1)(\Dhat\otimes\iota)(\Omega_\theta)=(\Dhat\otimes\iota)(\Omega_\theta)(\iota\otimes\Dhat)(F_\theta)$.
\end{proposition}

We don't have any conceptual explanation of this property. It is possible that this is true for any dual cocycle defined using a unitary quantization map as in~\cite[Section~2.2]{BGNT2}. In the present case the identity is not very difficult to verify by a direct computation as follows.

First of all, it is clearly enough to consider elementary K\"ahlerian Lie groups.

The next lemma is straightforward.

\begin{lemma}
For any $x\in\R$ and $\theta\in\R^*$, we have, with $g:=(x,0,0)\in G$:
$$
(\lambda_g\otimes\lambda_g)\Omega_\theta=\Omega_{e^{-2x}\theta}(\lambda_g\otimes\lambda_g).
$$
\end{lemma}

It follows that it suffices to prove the proposition for $\theta=\pm2$. In fact, by passing from $G$ to $G\rtimes_\alpha\Z/2\Z$, where $\alpha$ is the involutive automorphism defined at the Lie algebra level by $H\mapsto H$, $E\mapsto -E$, $X_i\mapsto X_{i+d}$ and $X_{i+d}\mapsto X_i$ for $i=1,\dots,d$, we can similarly relate $\Omega_2$ and $\Omega_{-2}$. Thus it is enough to consider $\theta=-2$.

From now on we write $\Omega$ and $F$ instead of $\Omega_{-2}$ and $F_{-2}$. The operator $\F_{P,-2}\colon L^2(G)\to L^2(G)$ from Appendix~\ref{ap:A} is up to a phase factor the partial Fourier transform $\F_P$ in $P$ coordinate:
$$
(\F_P f)(q,p)=(2\pi)^{-(d+1)/2}\int_{P} e^{-i\, p.p'}f(q,p')\,dp'.
$$
Therefore by Theorem~\ref{thm:A} we have $\Omega=(\F_{P}^*\otimes \F_{P}^*)U_\Phi(\F_{P}\otimes\F_{P})$.

Let us first consider the case $d=0$. Then
$$
\Phi(a_1,t_1;a_2,t_2)=\big(a_1-\tfrac12\arcsinh(e^{-2a_2}t_2),t_1;
a_2+\tfrac12\arcsinh(e^{-2a_1}t_1),t_2\big).
$$

\begin{lemma}\label{lem:O}
The map $\Phi$ is a diffeomorphism of $G\times G$ onto the proper subset
$$
\OO:=\Big\{(a_1,t_1;a_2,t_2)\mid \Big(1+\big(e^{-2a_1}t_1+e^{-2a_2}t_2\big)^2\Big)^{1/2}>e^{-2a_1}t_1-e^{-2a_2}t_2\Big\}\subset G\times G.
$$
\end{lemma}

\bp Consider the map $\phi\colon\R^2\to\R^2$, $\phi(a_1,a_2):=(\sinh(a_1)e^{-a_2},\sinh(a_2)e^{-a_1})$. It is easy to check that it is a diffeomorphism of $\R^2$ onto the set
$$
\big\{(a_1,a_2)\mid \big(1+(a_1-a_2)^2\big)^{1/2}>a_1+a_2\big\}\subset\R^2.
$$

Now, we have $\Phi(a_1,t_1;a_2,t_2)=(b_1,s_1,b_2,s_2)$ if and only if $t_i=s_i$ and
$$
\sinh(2b_1-2a_1)=-e^{-2a_2}s_2,\qquad \sinh(2b_2-2a_2)=e^{-2a_1}s_1,
$$
which in terms of the map $\phi$ is equivalent to
$$
\phi(2b_1-2a_1,2b_2-2a_2)=(-e^{-2b_2}s_2,e^{-2b_1}s_1).
$$
This gives the result.
\ep

It follows that $F=\Omega^*\Omega$ is indeed a proper projection, namely,
\begin{equation}\label{eq:k1}
F=(\F_{P}^*\otimes \F_{P}^*)M(\un_\OO)(\F_{P}\otimes\F_{P}),
\end{equation}
where $M(\un_\OO)$ is the operator of multiplication by the characteristic function $\un_\OO$ of $\OO$.

Next, recall that we have $\Dhat(x)=\hat W^*(1\otimes x)\hat W$, where the multiplicative unitary $\hat W\colon L^2(G\times G)\to L^2(G\times G)$ is defined by
$$
(\hat Wf)(g,h):=f(hg,h).
$$
The following lemma is again a straightforward computation.

\begin{lemma}
We have $\hat W=(\F_{P}^*\otimes \F_{P}^*)U_\alpha(\F_{P}\otimes\F_{P})$, where $U_\alpha f=f\circ\alpha$ and $\alpha\colon G\times G\to G\times G$ is the diffeomorphism defined by
$$
\alpha(a_1,t_1;a_2,t_2):=(a_1+a_2,t_1;a_2,t_2-e^{-2a_1}t_1).
$$
\end{lemma}

Note that the inverse of $\alpha$ is given by
$$
\alpha^{-1}(a_1,t_1;a_2,t_2)=(a_1-a_2,t_1;a_2,t_2+e^{2(a_2-a_1)}t_1).
$$

It follows that
\begin{align}
(\Dhat\otimes\iota)(\Omega)&=(\F_{P}^*\otimes\F_{P}^*\otimes \F_{P}^*)(U_\alpha)_{12}^*(U_{\Phi})_{23}(U_\alpha)_{12}(\F_{P}\otimes\F_{P}\otimes\F_{P}),\label{eq:k2}\\
(\iota\otimes\Dhat)(F)&=(\F_{P}^*\otimes\F_{P}^*\otimes \F_{P}^*)M(\un_{\alpha_{23}(\OO_{13})})(\F_{P}\otimes\F_{P}\otimes\F_{P}).\label{eq:k3}
\end{align}

From \eqref{eq:k1}--\eqref{eq:k3} we conclude that Proposition~\ref{prop:mystery} is equivalent (for $d=0$) to the following.

\begin{lemma}
We have
$$
\OO\times G=\{(g_1,g_2,g_3)\mid (\alpha_{12}\circ\Phi_{23}\circ\alpha_{12}^{-1})(g_1,g_2,g_3)\in \alpha_{23}(\OO_{13})\}.
$$
\end{lemma}

\bp
From Lemma~\ref{lem:O} we easily get that the points $(b_1,s_1;b_2,s_2;b_3,s_3)$ of $\alpha_{23}(\OO_{13})$ are characterized by the inequality
$$
\Big(1+\big(e^{-2b_1}s_1+e^{-2b_2}s_2+e^{-2b_3}s_3\big)^2\Big)^{1/2}>e^{-2b_1}s_1-e^{-2b_2}s_2-e^{-2b_3}s_3.
$$
Since
\begin{multline*}
(\alpha_{12}\circ\Phi_{23}\circ\alpha_{12}^{-1})(a_1,t_1;a_2,t_2;a_3,t_3)=\big(a_1-\tfrac12\arcsinh(e^{-2a_3}t_3),t_1;a_2-\tfrac12\arcsinh(e^{-2a_3}t_3),t_2;
\\ a_3+\tfrac12\arcsinh(e^{-2a_1}t_1+e^{-2a_2}t_2),t_3\big),
\end{multline*}
letting $x_i=e^{-2a_i}t_i$ we see that we have to show that the inequality
\begin{equation}\label{eq:ko1}
\Big(1+\big(x_1+x_2\big)^2\Big)^{1/2}>x_1-x_2
\end{equation}
is equivalent to
\begin{multline}\label{eq:ko2}
\Big(1+\big(e^{\arcsinh x_3}x_1+e^{\arcsinh x_3}x_2+e^{-\arcsinh(x_1+x_2)}x_3\big)^2\Big)^{1/2}\\>
e^{\arcsinh x_3}x_1-e^{\arcsinh x_3}x_2-e^{-\arcsinh(x_1+x_2)}x_3.
\end{multline}

Put $\displaystyle a:=e^{\arcsinh(x_1+x_2)}>0$ and $b:=x_1-x_2$. Then $x_1+x_2=(a-a^{-1})/2$, so~\eqref{eq:ko1} becomes
\begin{equation}\label{eq:ko1a}
\frac{a+a^{-1}}{2}>b.
\end{equation}
Similarly, letting $c:=e^{\arcsinh x_3}>0$, we get that~\eqref{eq:ko2} becomes
$$
\frac{ac+a^{-1}c^{-1}}{2}>bc-\frac{a^{-1}c-a^{-1}c^{-1}}{2},
$$
which is obviously equivalent to~\eqref{eq:ko1a} for all $a>0$, $b\in\R$ and $c>0$.
\ep

This completes the proof of Proposition~\ref{prop:mystery} for $d=0$. The case $d>0$ can be easily deduced from this. Indeed, from the definition of $\Phi$ we see that $\tilde\OO:=\Phi(G\times G)$ is the preimage of the set $\OO$ under the map~$\pi$ defined by $\pi(a_1,v_1,t_1;a_2,v_2,t_2):=(a_1,t_1;a_2,t_2)$. It is also easy to check that the transformation $\tilde\alpha\colon G\times G\to G\times G$ defining the unitary $(\F_{P}\otimes \F_{P})\hat W(\F_{P}^*\otimes\F_{P}^*)$ satisfies $\pi\circ\tilde\alpha=\alpha\circ\pi$. From this we see that the coordinates $v$ do not play any role, so the result follows from the case $d=0$.

\medskip

Therefore $\Omega_\theta$ satisfies assumptions (1)--(3) from Section~\ref{sec:def}, so we get a well-behaved notion of $\Omega_\theta$-deformation. The main result of the paper - Theorem~\ref{main-stuff2} - then holds as stated, and essentially the only change in the proof is that the operator $\hat W_{\Omega_\theta}\Omega_\theta$ gets replaced by the operator $(J\otimes \hat J)\Omega_\theta\hat W^*(J\otimes\hat J)$ everywhere.

In more detail, Corollary~\ref{cregular} and Proposition~\ref{pmodular} and the discussions preceding them are no longer meaningful and should be removed. The rest of the paper remains unchanged, with identity~\eqref{emult} understood as a definition of the operator~$\hat W_{\Omega_\theta}\Omega_\theta$. Note that the place where we need special properties of $\Omega_\theta$ is the first identity in the proof of Lemma~\ref{l:I-have-no-money}, which now holds by Proposition~\ref{prop:eq3.1}.

\bigskip


\bigskip

\begin{thebibliography}{99}

\bibitem{BC}
S. Baaj and J. Crespo, {\em Equivalence monoidale de groupes quantiques et K-th\'eorie bivariante}, preprint arXiv: 1507.06808v1 [math.OA].

\bibitem{Bi}
P. Bieliavsky, {\em Strict quantization of solvable symmetric spaces},
J. Symplectic Geom. {\bf 1} (2002), no.~2, 269--320.

\bibitem{BG}
P. Bieliavsky and V. Gayral, {\em Deformation Quantization for Actions of K\"ahlerian Lie Groups},
 Mem. Amer. Math. Soc. {\bf 236} (2015),   no.~1115.

\bibitem{BGNT2}
P. Bieliavsky, V. Gayral, S. Neshveyev and L. Tuset, {\em Quantization of subgroups of the affine group}, preprint arXiv:1906.01889v2 [math.OA].

\bibitem{BNS}
J. Bhowmick, S. Neshveyev and A. Sangha, {\em Deformation of operator algebras by Borel cocycles}, J. Funct. Anal. {\bf 265} (2013), no.~6, 983--1001.

\bibitem{DC}
K. De Commer, {\em Galois objects and cocycle twisting for locally compact quantum groups}, J. Operator Theory {\bf 66}  (2011),  no. 1, 59--106.

\bibitem{Kas}
P. Kasprzak, {\em Rieffel deformation via crossed products}, J. Funct. Anal. {\bf 257}  (2009),  no.~5, 1288--1332.

\bibitem{N}
S. Neshveyev, {\em Smooth crossed products of Rieffel's deformations},
Lett. Math. Phys. {\bf 104} (2014), no.~3, 361--371.

\bibitem{NTdef2}
S. Neshveyev and L. Tuset, {\em Deformation of C$^*$-algebras by cocycles on locally compact quantum groups}, Adv. Math. {\bf 254} (2014), 454--496.

\bibitem{PS}
I. {Pyateskii-Shapiro},
{\it Automorphic functions and the geometry of classical domains},
\newblock Translated from the Russian. Mathematics and Its Applications, Vol.
 {\bf  8}, Gordon and Breach Science Publishers, New York (1969).

\bibitem{Ri1}
M.A. Rieffel, {\em Deformation quantization for actions of ${\mathbb R}^d$}, Mem. Amer. Math. Soc. {\bf 106}  (1993),  no.~506.

\end{thebibliography}
\end{document}